\documentclass[12 pt]{article}

\usepackage[cp1251]{inputenc}
\usepackage[french, english]{babel}
\usepackage{morefloats}
\usepackage{amssymb} 
\usepackage{amsfonts} 
\usepackage{amsmath}
\usepackage{mathrsfs}
\usepackage{amsthm} 
\usepackage{epsfig} 
\usepackage{color}
\usepackage{amscd}
\usepackage[all]{xy}
\usepackage{subfigure}
\usepackage{bookmark}
\newtheorem{lemma}{Lemma}[section]
\newtheorem{teo}[lemma]{Theorem}
\newtheorem{prop}[lemma]{Proposition}
\newtheorem{cor}[lemma]{Corollary}

\theoremstyle{definition}
\newtheorem{defn}[lemma]{Definition}

\newtheorem{quest}[lemma]{Question}

\theoremstyle{remark}
\newtheorem{rem}[lemma]{Remark}

\newcommand{\calR} {\ensuremath {\mathcal{R}}}

\newcommand{\calM} {\ensuremath {\mathcal{M}}}

\newcommand{\calB} {\ensuremath {\mathcal{B}}}

\newcommand{\calT} {\ensuremath {\mathcal{T}}}

\newcommand{\prf}[1]{\vspace{2pt}\noindent\textit{Proof of #1}.\ }
\newcommand{\prfend}{{\hfill\hbox{$\square$}\vspace{2pt}}}

\usepackage[letterpaper,left=0.9in,right=0.9in,top=1in,bottom=1in]{geometry}
\usepackage{mathtools}

\DeclarePairedDelimiter\floor{\lfloor}{\rfloor}

\begin{document}
\title{Symmetries of hyperbolic 4-manifolds}
\author{Alexander Kolpakov \& Leone Slavich}
\date{}
\maketitle

\selectlanguage{french} 
\begin{abstract}
Pour chaque groupe $G$ fini, nous construisons des premiers exemples explicites de $4$-vari\'et\'es non-compactes compl\`etes arithm\'etiques hyperboliques $M$, \`a volume fini, telles que $\mathrm{Isom}\, M \cong G$, ou $\mathrm{Isom}^{+}\, M \cong G$. Pour y parvenir, nous utilisons essentiellement la g\'eom\'etrie de poly\`edres de Coxeter dans l'espace hyperbolique en dimension quatre, et aussi la combinatoire de complexes simpliciaux. 

\c{C}a nous permet d'obtenir une borne sup\'{e}rieure universelle pour le volume minimal d'une $4$-vari\'{e}t\'{e} hyperbolique ayant le groupe $G$ comme son groupe d'isom\'etries, par rapport de l'ordre du groupe. Nous obtenons aussi des bornes asymptotiques pour le taux de croissance, par rapport du volume, du nombre de $4$-vari\'{e}t\'{e}s hyperboliques ayant $G$ comme le groupe d'isom\'etries.
\end{abstract}

\selectlanguage{english} 
\begin{abstract}
In this paper, for each finite group $G$, we construct the first explicit examples of non-compact complete finite-volume arithmetic hyperbolic $4$-manifolds $M$ such that $\mathrm{Isom}\,M \cong G$, or $\mathrm{Isom}^{+}\,M \cong G$. In order to do so, we use essentially the geometry of Coxeter polytopes in the hyperbolic $4$-space, on one hand, and the combinatorics of simplicial complexes, on the other. 

This allows us to obtain a universal upper bound on the minimal volume of a hyperbolic $4$-manifold realising a given finite group $G$ as its isometry group in terms of the order of the group. We also obtain asymptotic bounds for the growth rate, with respect to volume, of the number of hyperbolic $4$-manifolds having a finite group $G$ as their isometry group.
\end{abstract}

\tableofcontents

\section{Introduction}

In this paper we give the first explicit examples of complete hyperbolic manifolds with given isometry group in dimension four. All our manifolds have finite volume and are arithmetic, by construction. Our interest in constructing explicit and feasible examples is motivated by the work of M. Belolipetsky and A. Lubotzky \cite{BL}, which shows that for any finite group $G$ and any dimension $n\geq 2$, there exists a complete, finite volume, hyperbolic non-arithmetic manifold $M$\footnote{in fact, \cite{BL} shows that there are infinitely many manifolds $M$ with $\mathrm{Isom}\,M \cong G$.}, such that $\mathrm{Isom}\,M \cong G$. This statement was proved earlier, with various methods, for $n=2$ in \cite{E} and \cite{G}, for $n=3$ first in \cite{Kojima}, and then, in a more general context, in \cite{FM}. The case of the trivial group $G=\{e\}$ was considered in \cite{LR}. 

The construction of such a manifold $M$ in \cite{BL} utilises the features of arithmetic group theory, similar to the preceding work by D. Long and A. Reid \cite{LR}, and the subgroup growth theory, which provides a probabilistic argument in proving the existence of $M$. 

In the present paper we use the methods of Coxeter group theory and combinatorics of simplicial complexes, close to the techniques of \cite{FM} and \cite{KM}. These methods allow us to construct manifolds with highly controllable geometry and we are able to estimate their volume in terms of the order of the group $G$. The main results of the paper read as follows:

\begin{teo}
Given a finite group $G$ there exists an arithmetic, non-orientable, four-di\-men\-sio\-nal, complete, finite-volume, hyperbolic manifold $M$, such that $\mathrm{Isom}\,M \cong G$.
\end{teo}

\begin{teo}
Given a finite group $G$ there exists an arithmetic, orientable, four-di\-men\-sio\-nal, complete, finite-volume, hyperbolic manifold $M$, such that $\mathrm{Isom}^{+}\,M \cong G$.
\end{teo}

As a by-product of our construction we obtain that

\begin{teo}
The group $G$ acts on the manifold $M$ freely.
\end{teo}

We also give an upper bound on the volume of the manifold in terms of the order of the group $G$, giving a partial answer to a question first asked in \cite{BL}:

\begin{teo}
Let the group $G$ have rank $m$ and order $n$. Then in the above theorems we have $\mathrm{Vol}\,M \leq C\cdot n\cdot m^2$, where the constant $C$ does not depend on $G$.
\end{teo}

The paper is organised as follows: first we discuss the initial ``building block'' of our construction, which comes from assembling six copies of the ideal hyperbolic rectified $5$-cell, and prove that this object is combinatorially equivalent to the standard $4$-dimensional simplex. Then, given a $4$-dimensional simplicial complex $\mathcal{T}$, called a triangulation, we associate a non-orientable manifold $M_{\mathcal{T}}$ with it. We also prove that our manifolds $M_\mathcal{T}$, up to an isometry, are in a one-to-one correspondence with the set of triangulations, up to a certain combinatorial equivalence. 

We show how the structure of the triangulation $\mathcal{T}$ encodes the geometry and topology of the manifold $M_\mathcal{T}$: the maximal cusp section of $M_\mathcal{T}$ is uniquely determined by $\mathcal{T}$, as well as the isometry group $\mathrm{Isom}\,M_\mathcal{T}$. 

Finally, we construct a triangulation $\mathcal{T}$ with a given group $G$ of combinatorial automorphisms, and thus obtain the desired manifold $M := M_\mathcal{T}$ with $\mathrm{Isom}\,M \cong G$. Its orientable double cover produces a manifold $\widetilde{M}$, such that $\mathrm{Isom}^+\,\widetilde{M} \cong G$. Finally, we estimate the volume of $M$, which is a direct consequence of our construction. 

In the last section, we show that the number of manifolds having a given finite group $G$ as their isometry group, grows super-exponentially with respect to volume. More precisely, defining $\rho_G(V) = \#\{\, M\, |\, \mathrm{Isom}\,M \cong G \mbox{ and } \mathrm{Vol}\,M \leq V \}$, we prove:

\begin{teo}\label{teo:superexponential1}
For any finite group $G$, there exists a $V_0 > 0$ sufficiently large such that for all $V \geq V_0$ we have
$\rho_G(V) \geq C^{\, V \log V}$, for some $C>1$ independent of $G$.
\end{teo}
To do so, we combine our construction with some counting results \cite{Bo, KSV} on the number of trivalent graphs on $n$ vertices. 

\medskip

\textbf{Acknowledgements.} The authors are grateful to FIRB 2010 "Low-dimensional geometry and topology" and the organisers of the workshop ``Teichm\"uller theory and surfaces in 3-manifolds'' during which the most part of the paper was written. The authors received financial support from FIRB (L.S., FIRB project no. RBFR10GHHH-003) and the Swiss National Science Foundation (A.K., SNSF project no. P300P2-151316). Also, the authors are grateful to Bruno Martelli (Universit\`a di Pisa), Ruth Kellerhals (Universit\'e de Fribourg), Marston Conder (University of Auckland), Sadayoshi Kojima (Tokyo Institute of Technology), Makoto Sakuma (Hiroshima University) and Misha Belolipetsky (IMPA, Rio de Janeiro) for fruitful discussions and useful references. A.K. is grateful to Waseda University and, personally, to Jun Murakami for hospitality during his visit in autumn 2014, when a part of this work was finished. 

\section{The rectified 5-cell}

Below, we describe the main building ingredient of our construction, \textit{the rectified $5$-cell}, which can be realised as a non-compact finite-volume hyperbolic $4$-polytope. First, we start from its Euclidean counterpart, which shares the same combinatorial properties.  

\begin{defn}\label{defn:5cell}
The \emph{Euclidean rectified $5$-cell} $\mathcal{R}$ is the convex hull in $\mathbb{R}^5$ of the set of $10$ points whose coordinates are obtained as all possible permutations of those of the point $(1,1,1,0,0)$.
\end{defn}

The rectified $5$-cell has ten facets ($3$-dimensional faces) in total. Five of these are regular octahedra. They lie in the affine planes defined by the equations 
\begin{equation}
\sum_{i=1}^5 x_i=3,\; x_j=1, \mbox{ for each } j\in \{1,2,3,4,5\},
\end{equation}
and are naturally labelled by the number $j$. 

The other five facets are regular tetrahedra. They lie in the affine hyperplanes given by the equations 
\begin{equation}
\sum_{i=1}^5 x_i=3,\; x_j=0, \mbox{ for each } j\in \{1,2,3,4,5\},
\end{equation}
and are also labelled by the number $j$.

Also, the polytope $\mathcal{R}$ has $30$ two-dimensional triangular faces, $30$ edges and $10$ vertices.

We note the following facts about the combinatorial structure of $\mathcal{R}$:
\begin{enumerate}
\item each octahedral facet $F$ has a red/blue chequerboard colouring, such that $F$ is adjacent to any other octahedral facet along a red face, and to a tetrahedral facet along a blue face;
\item a tetrahedral facet having label $j\in \{1,2,3,4,5\}$ is adjacent along its faces to the four octahedra with labels $k\in \{1,2,3,4,5\}$, with $k$ different from $j$;
\item the tetrahedral facets meet only at vertices and their vertices comprise all those of $\mathcal{R}$.
\end{enumerate}

\begin{rem} Another way to construct the rectified $5$-cell is to start with a regular Euclidean $4$-dimensional simplex $S_4$ and take the convex hull of the midpoints of its edges. This is equivalent to truncating the vertices of $S_4$, and enlarging the truncated regions until they become pairwise tangent along the edges of $S_4$.

With this construction, it is easy to see that the symmetry group of $\mathcal{R}$ is isomorphic to the symmetry group of $S_4$, which is known to be $\mathfrak{S}_5$, the group of permutations of a set of five elements.
\end{rem}

\begin{defn}\label{defn:hyperbolic-5-cell}
Like any other uniform Euclidean polytope, the rectified $5$-cell has a hyperbolic ideal realisation, which may be obtained in the following way:
\begin{enumerate}
\item normalise the coordinates of the vertices of $\mathcal{R}$ so that they lie on the unit sphere $\mathbb{S}^3\subset \mathbb{R}^4$;
\item interpret $\mathbb{S}^3$ as the boundary at infinity of the hyperbolic $4$-space $\mathbb{H}^4$ in the Klein-Beltrami model.
\end{enumerate}
The convex hull of the vertices of $\mathcal{R}$ now defines an ideal polytope in $\mathbb{H}^4$, that we call the \emph{ideal hyperbolic rectified $5$-cell}.
\end{defn}

With a slight abuse of notation, we continue to denote the ideal hyperbolic rectified $5$-cell by $\mathcal{R}$.

\begin{rem}
The vertex figure of the ideal hyperbolic rectified $5$-cell is a right Euclidean prism over an equilateral triangle, with all edges of equal length. At each vertex, there are three octahedra meeting side-by-side, corresponding to the square faces, and two tetrahedra, corresponding to the triangular faces.

The dihedral angle between two octahedral facets is therefore equal to $\pi/3$, while the dihedral angle between a tetrahedral and an octahedral facet is equal to $\pi/2$. 
\end{rem}

\begin{rem}
The volume $v_{\mathcal{R}}$ of the rectified $5$-cell equals $2\pi^2/9$, as computed in Appendix A.
\end{rem}

\section{The building block}\label{sec:buildblock}

In this section, we produce a building block $\mathcal{B}$, which is the second stage of our construction. We show that $\mathcal{B}$ is in fact a non-compact finite-volume hyperbolic manifold with totally geodesic boundary, and then study its isometry group $\mathrm{Isom}\,\mathcal{B}$.

\begin{figure}[ht]
\centering
\includegraphics[width=0.5\textwidth]{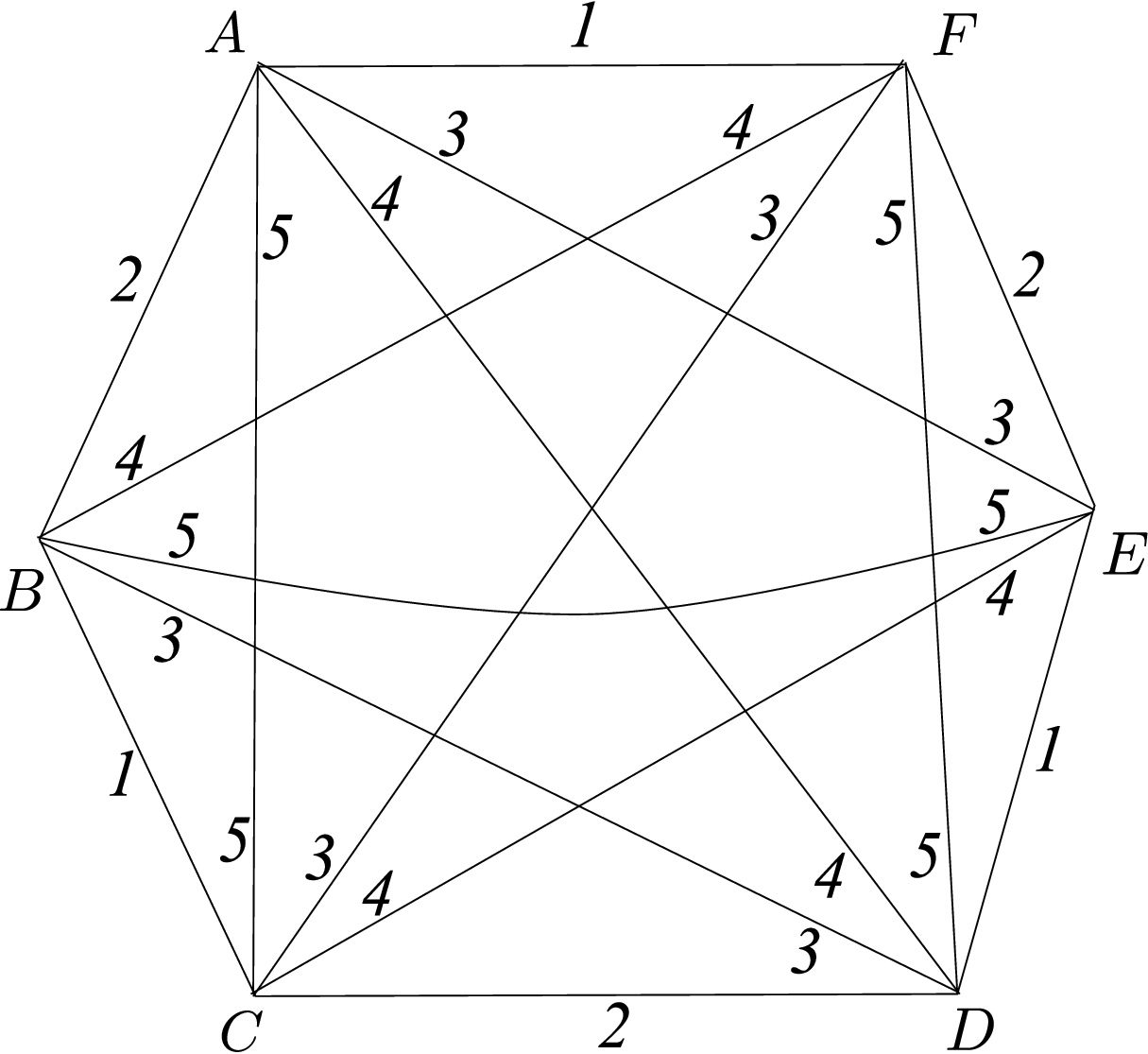}
\caption{The glueing graph $\Gamma$ for the building block $\mathcal{B}$.}\label{fig:GlueingGraphBlock}
\end{figure}

Let us consider six copies of the ideal hyperbolic rectified $5$-cell $\mathcal{R}$, labelled by the letters $A$, $B$, $C$, $D$, $E$, $F$. Recall that each of the octahedral facets of these copies of $\mathcal{R}$ is naturally labelled by an integer $i\in\{1,2,3,4,5\}$. Let us pair all the octahedral facets according to the glueing graph $\Gamma$ in Fig.~\ref{fig:GlueingGraphBlock}, always using the identity as a pairing map. A label on each edge specifies which octahedral facets are paired together.

\begin{rem}
The edge labelling of the graph $\Gamma$ is, up to a permutation of the numbers, the only possible one with five numbers on the complete graph $K_6$ on six vertices.
\end{rem}

Let us denote by $\mathcal{B}$ the resulting object. The following proposition clarifies its nature.

\begin{prop}
The building block $\mathcal{B}$ is a complete, non-orientable, finite-volume, hyperbolic manifold with totally geodesic boundary. Its volume equals
\begin{equation}
v_{\mathcal{B}}=6\cdot v_{\mathcal{R}}=4\pi^2/3.
\end{equation}
\end{prop}

\begin{proof}
We need to check that the natural hyperbolic structures on the six copies of $\mathcal{R}$ match together under the glueing to give a complete hyperbolic structure on the whole manifold $\mathcal{B}$. In order to do so, it suffices to check that the pairing maps glueing together the respective Euclidean vertex figures along their square faces actually produce Euclidean $3$-manifolds (with totally geodesic boundary) as cusp sections.

Indeed, the pairing maps that define $\mathcal{B}$ produce ten cusps, which are in a one-to-one correspondence with the vertices of $\mathcal{R}$. The cusp section is a trivial $I$-bundle over a flat Klein bottle $K$, tessellated by six equilateral triangles, each one coming from its own copy of $\mathcal{R}$, as shown in Fig.~\ref{fig:cusps}.

\begin{figure}[ht]
\centering
\includegraphics[width=0.35\textwidth]{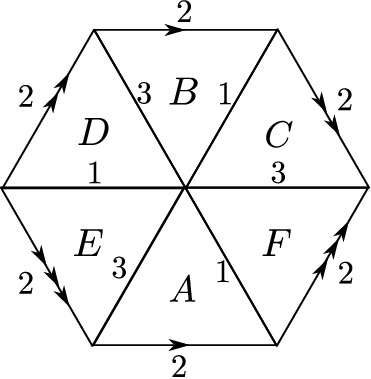}
\caption{The cusp section of $\mathcal{B}$ corresponding to the vertex $(1,1,1,0,0)$ of $\mathcal{R}$ is a trivial $I$-bundle over the Klein bottle shown above.}\label{fig:cusps}
\end{figure}

The volume of the building block $\mathcal{B}$ is equal to $6\cdot v_{\mathcal{R}}$, since $\mathcal{B}$ is built by glueing together six copies of $\mathcal{R}$.
\end{proof}

\begin{prop}
The block $\mathcal{B}$ has five totally geodesic boundary components, all isometric to each other.
\end{prop}

\begin{proof}
The boundary components are in a one-to-one correspondence with the tetrahedral facets of the rectified $5$-cell $\mathcal{R}$, and are therefore naturally labelled by an integer $i\in\{1,2,3,4,5\}$. The glueing graph $\hat{\Gamma}$ for the boundary component labelled by $i$ is obtained from the glueing graph $\Gamma$ in Fig.~\ref{fig:GlueingGraphBlock}, by removing all the edges labelled by $i$, as shown in Fig.~\ref{fig:boundarygraph}. In this case, each vertex corresponds to a copy of the regular ideal hyperbolic tetrahedron and the pairing maps are once more induced by the identity.

\begin{figure}[ht]
\centering
\includegraphics[width=0.3\textwidth]{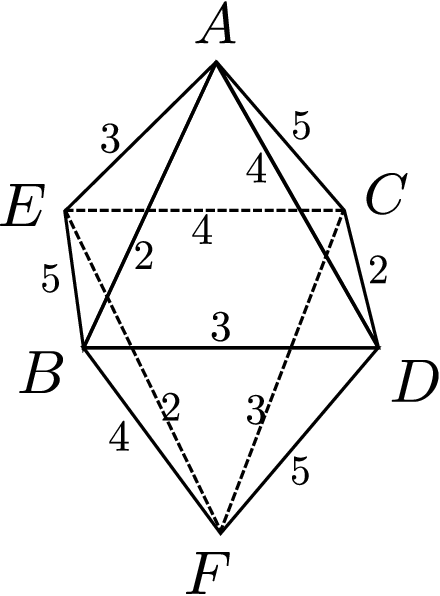}
\caption{The glueing graph $\hat{\Gamma}$ for the boundary component of $\mathcal{B}$ labelled $1$, together with the edge labels.}\label{fig:boundarygraph}
\end{figure}

The resulting labelled graphs are all isomorphic, and therefore all boundary components are isometric.
\end{proof}

From now on, we will denote by $\mathcal{D}$ the hyperbolic $3$-manifold isometric to the boundary components of $\mathcal{B}$, and by $\widetilde{\mathcal{D}}$ its orientable double cover.

\begin{rem}
The orientable double cover $\widetilde{\mathcal{D}}$ of a boundary component of $\mathcal{B}$ is the complement of a link with four components. The link is depicted in \cite[p. 148]{ASh}, at the entry $n=4$, $\sigma(n)=6$.
\end{rem}

\begin{defn}\label{defn:bounding}
Let $\mathcal{M}$ be an orientable, complete, finite-volume hyperbolic $n$-manifold without boundary. We say that $\mathcal{M}$ \textit{geometrically bounds} if there exists an orientable, complete, finite-volume, hyperbolic $(n+1)$-manifold $\mathcal{X}$ such that:
\begin{enumerate}
\item $\mathcal{X}$ has only one totally geodesic boundary component;
\item the boundary $\partial \mathcal{X}$ is isometric to $\mathcal{M}$.
\end{enumerate} 
\end{defn}

\begin{prop}\label{prop:linkbounds}
The hyperbolic $3$-manifold $\widetilde{\mathcal{D}}$ bounds geometrically.
\end{prop}

\begin{proof}
The orientable double cover $\widetilde{\mathcal{B}}$ of the block $\mathcal{B}$ has five totally geodesic boundary components, all isometric to $\widetilde{\mathcal{D}}$. Indeed, both $\mathcal{B}$ and its boundary components are non-orientable, so each connected component of $\partial \mathcal{B}$ lifts to a single connected component of $\partial \widetilde{\mathcal{B}}$. 

We can identify four of the boundary components of $\widetilde{\mathcal{B}}$ in pairs, using \emph{any} orientation-reversing isometry as the respective glueing map. Indeed, $\widetilde{\mathcal{D}}$ has an orientation-reversing involution corresponding to the orientable double-cover. Thus, we obtain an orientable complete finite-volume hyperbolic $4$-manifold with a single totally geodesic boundary component isometric to $\widetilde{\mathcal{D}}$.
\end{proof}

\begin{rem}
The hyperbolic volume of $\widetilde{\mathcal{D}}$ is $\approx 12.1792\dots$. Therefore $\widetilde{\mathcal{D}}$ has the lowest volume amongst known geometrically bounding hyperbolic $3$-manifolds, cf. \cite{RT}, \cite{S1}. 
\end{rem}

\begin{rem}
The ratio between the volume of the building block $\calB$ and the volume of its boundary $\partial \calB$ is $$\frac{\text{Vol}\, \calB}{\text{Vol}\, \partial \calB}= 0.43219...$$
 
Moreover it is clear that the same ratio holds for the orientable double cover $\widetilde{\calB}$ of the building block. These manifolds are therefore the first explicit $4$-dimensional examples realising Miyamoto's lower bound \cite[Theorem~4.2, Proposition~4.3]{M}, which relates the volume of a hyperbolic manifold with the volume of its totally geodesic boundary.
\end{rem}

\subsection{Combinatorial equivalence}

We shall establish a \emph{combinatorial equivalence} between the building block $\mathcal{B}$ and the standard $4$-dimensional simplex $S_4$. In particular, the following one-to-one correspondences hold:
\begin{enumerate}
\item $\{\mbox{Boundary components of }\mathcal{B}\}\leftrightarrow\{\mbox{Facets of }S_4\}$
\item $\{\mbox{Cusps of }\mathcal{B}\}\leftrightarrow\{\mbox{2-dimensional faces of }S_4\}$.
\end{enumerate}
In order to see that these correspondences hold, it is enough to notice that the block $\mathcal{B}$ can be decomposed into six copies of the ideal hyperbolic rectified $5$-cell $\mathcal{R}$, and that any of these copies will have its vertices in a one-to-one correspondence with the cusps of $\mathcal{B}$, and its tetrahedral facets in a one-to-one correspondence with the boundary components of $\mathcal{B}$. Moreover, viewing the $5$-cell $\mathcal{R}$ as the result of truncation of a $4$-dimensional simplex $S_4$, the vertices of $\mathcal{R}$ naturally correspond to the edges of $S_4$ and the tetrahedral facets naturally correspond to the vertices of $S_4$. By considering the dual polytope to $S_4$ (which is again $S_4$), we obtain all the desired correspondences.

The above combinatorial equivalence between the strata of the simplex $S_4$ and the geometric compounds of the block $\mathcal{B}$ allows us to describe the isometry group of $\mathcal{B}$ and that of its boundary component $\mathcal{D}$.

\begin{prop}\label{prop:isometry-boundary}
There is an isomorphism between the group $\mathrm{Isom}\,\mathcal{D}$ of isometries of the boundary manifold $\mathcal{D}$ and the group $\mathfrak{S}_4$ of symmetries of a tetrahedron.
\end{prop}

\begin{proof}
We begin by showing that, for any permutation $\sigma\in \mathfrak{S}_4$ of the edge colours of $\hat{\Gamma}$ (which are $2,3,4,5$, as in Figure \ref{fig:boundarygraph}), there is a unique automorphism $\phi_{\sigma}$ of $\hat{\Gamma}$ (viewed as an unlabelled graph) which permutes the labels on the edges in the way defined by $\sigma$. Without loss of generality, we can suppose that $\sigma$ is a transposition of two labels, for instance the transposition of the labels $3$ and $4$ on the graph $\hat{\Gamma}$ in Fig.~\ref{fig:boundarygraph}. The automorphism $\phi_{\sigma}$ is defined by the following map of the vertices:
\begin{center}
$A\leftrightarrow F$, 

$E\leftrightarrow B$, 

$D\leftrightarrow C$.
\end{center}

The uniqueness of $\phi_{\sigma}$ follows from the fact that the group of automorphisms of $\hat{\Gamma}$ \emph{as a labelled graph} is trivial. To see this, notice that the vertices and edges of $\hat{\Gamma}$  form the one-skeleton of an octahedron. Any automorphism $\phi$ of $\hat{\Gamma}$ as a labelled graph is required to fix all pairs of opposite faces, since $\phi$ preserves cycles of vertices of length $3$, which correspond to the faces of the octahedron, and only opposite faces share the same labels on their edges.

Let us suppose that two opposite faces $F_1$ and $F_2$ are fixed by an automorphism $\phi$ (in the sense that $\phi(F_1)=
F_1$ and $\phi(F_2)=F_2$). Then $\phi$ is necessarily the identity. If instead $\phi(F_1)=F_2$ and $\phi(F_2)=F_1$, the image of each vertex under $\phi$ is uniquely determined by the respective edge labels, but there will always be a couple of vertices $v$ and $w$ in $\hat{\Gamma}$ that share an edge, and whose images under $\phi$ are a couple of opposite vertices of the octahedron. Therefore, $\phi(v)$ and $\phi(w)$ do not share an edge, which is a contradiction.

Now we notice that every vertex $v$ of $\hat{\Gamma}$ corresponds to a tetrahedron $T_v$ that tessellates $\mathcal{D}$, and that every edge of $\hat{\Gamma}$ adjacent to $v$ corresponds to a unique triangular face of $T_v$. Given $\sigma\in \mathfrak{S}_4$, we map each tetrahedron $T_v$ to $T_{\phi_{\sigma}(v)}$, respecting the pairing on the triangular facets defined by the permutation $\sigma$. This defines an injective homomorphism from $\mathfrak{S}_4$ to $\mathrm{Isom}\,\mathcal{D}$. Also, this homomorphism is surjective, which follows from the fact that every isometry of $\mathcal{D}$ has to fix its Epstein-Penner decomposition \cite{EP} into regular ideal hyperbolic tetrahedra. 
\end{proof}

\begin{prop}\label{prop:isometry-block}
There is an isomorphism between the group $\mathrm{Isom}\,\mathcal{B}$ of isometries of the building block $\mathcal{B}$ and the group $\mathfrak{S}_5$ of symmetries of a $4$-dimensional simplex $S_4$.
\end{prop}

\begin{proof}
An isometry of $\mathcal{B}$ acts on the set of its five boundary components as a (possibly, trivial) permutation. Thus, we have a natural homomorphism from $\mathrm{Isom}\, \mathcal{B}$ to $\mathfrak{S}_5$. This homomorphism is in fact an isomorphism. As a first step, we notice that any isometry of the building block $\cal{B}$ has to preserve its decomposition into copies of the rectified $5$-cell $\cal{R}$. We postpone the proof, since this is a particular case (Corollary \ref{cor:buildingblock}) of a much more general statement (Proposition \ref{prop:Epstein-Penner}) which we will prove later on.  

Because of this fact, any isometry of the block $\cal{B}$ has to induce an automorphism of the glueing graph $\hat{\Gamma}$, which perhaps permutes the edge labels.
Indeed, we have the inclusion $\hat{\Gamma}\subset \Gamma$ and we know that every transposition of the edge colours of $\hat{\Gamma}$ is obtained by a unique automorphism of $\hat{\Gamma}$ as an unlabelled graph, as shown in the proof of Proposition~\ref{prop:isometry-boundary}. This automorphism extends to a \emph{unique} automorphism of the whole graph $\Gamma$, which necessarily preserves one of the labels on the edges. Moreover, every automorphism of $\Gamma$ as a labelled graph induces an automorphism of $\hat{\Gamma}$ as a labelled graph, and therefore has to be the identity.
\end{proof}

\subsection{The maximal cusp section}
The ideal hyperbolic rectified $5$-cell $\mathcal{R}$ has a canonical maximal cusp section. It is obtained by placing the vertices of the Euclidean rectified $5$-cell on the boundary at infinity of $\mathbb{H}^4$, as in Definition~\ref{defn:hyperbolic-5-cell}, and expanding uniformly (with respect to the Euclidean metric) the horospheres centred at the vertices until they all become pairwise tangent. 

With this choice, the edges of the Euclidean vertex figure all have length one. To see this, notice that the intersection of the horospheres constructed above with any $2$-face of $\calR$ is given by three horocyclic segments in an ideal hyperbolic triangle. These segments intersect pairwise at their endpoints on the edges of the triangle. In each ideal triangle, there is only one such collection of segments, and they all have length one. Also, we observe that exactly these segments form the edges of the equilateral triangles and squares that constitute the faces in each vertex figure of the rectified $5$-cell. Thus, the edges of the maximal cusp section necessarily have length one.

When we build the block $\mathcal{B}$ by glueing together six copies of $\mathcal{R}$, the maximal cusp sections of each copy are identified isometrically along their square faces in order to produce the maximal cusp section of $\mathcal{B}$.

\section{Hyperbolic 4-manifolds from triangulations}
Below, we produced a hyperbolic $4$-manifold $\mathcal{M}$ from the combinatorial data carried by a $4$-dimensional triangulation $\mathcal{T}$ and describe its isometry group $\mathrm{Isom}\,\mathcal{M}$.

\begin{defn}\label{defn:triangulation}
A \textit{$4$-dimensional triangulation} $\mathcal{T}$ is a pair 
\begin{equation}
(\{\Delta_i\}_{i=1}^{2n}, \{g_j\}_{j=1}^{5n}),
\end{equation} 
where $n$ is a positive natural number, the $\Delta_i$'s are copies of the standard $4$-dimensional simplex $S_4$, and the $g_j$'s are a complete set of simplicial pairings between the $10n$ facets of all $\Delta_i$'s.
\end{defn}

\begin{defn}\label{defn:orientable-triangulation}
A triangulation is \textit{orientable} if it is possible to choose an orientation for each tetrahedron $\Delta_i$, $i=1,\dots,2n$, so that all pairing maps between the facets are orientation-reversing.
\end{defn}

\begin{defn}\label{defn:combinatorial-equivalence}
A \textit{combinatorial equivalence} between two $4$-dimensional triangulations $\mathcal{T}=(\{\Delta_i\}_{i=1}^{2n}, \{g_j\}_{j=1}^{5n})$ and $\mathcal{T}'=(\{\Delta'_i\}_{i=1}^{2n}, \{g'_j\}_{j=1}^{5n})$ is a set of simplicial maps $\phi_{kl}: \Delta_k\rightarrow \Delta'_l$ which induces a one-to-one correspondence between the pairings of $\mathcal{T}$ and the pairings of $\mathcal{T}'$.

Given a triangulation $\mathcal{T}$, the group of combinatorial equivalences of $\mathcal{T}$ is denoted by $\mathrm{Aut}\, \mathcal{T}$, and an element of such group is called an \emph{automorphism} of $\mathcal{T}$.
\end{defn}

By virtue of the combinatorial equivalence between the block $\mathcal{B}$ and the $4$-di\-men\-sio\-nal simplex $S_4$ deduced in the previous section, we can encode the pairings between the boundary components of several copies of $\mathcal{B}$ by using simplicial face pairings between the facets of copies of $S_4$. This allows us to produce a hyperbolic $4$-manifold from the data carried by a $4$-dimensional triangulation. 

The construction is as follows:
\begin{enumerate}
\item given a triangulation $\mathcal{T}=(\{\Delta_i\}_{i=1}^{2n}, \{g_j\}_{j=1}^{5n})$, associate with each $\Delta_i$ a copy $\mathcal{B}_i$ of the building block;
\item a face pairing $g_{kl}$ between the facets $F$ and $G$ of the simplices $\Delta_k$ and $\Delta_l$ defines a \textit{unique} isometry between the respective boundary components $\mathcal{D}_F$ and $\mathcal{D}_G$ of the blocks $\mathcal{B}_k$ and $\mathcal{B}_l$, as in Proposition~\ref{prop:isometry-boundary};
\item identify all boundary components of the blocks $\mathcal{B}_i$, $i=1\dots,2n$ using the isometries defined by the pairings $g_j$, $j=1\dots 5n$, to produce $M_{\mathcal{T}}$.
\end{enumerate}

The nature of the above constructed object $M_{\mathcal{T}}$ is clarified by the following proposition. 

\begin{prop}\label{prop:hyperbolic-structure}
If $\mathcal{T}$ is a $4$-dimensional triangulation, then $M_{\mathcal{T}}$ is a non-orientable, non-compact, complete, finite-volume, arithmetic hyperbolic $4$-manifold.
\end{prop}

\begin{proof}
Since the copies of the block $\mathcal{B}$ are glued together via isometries of their totally geodesic boundary components, their hyperbolic structures match together to give a hyperbolic structure on $M_{\mathcal{T}}$. Its volume equals $2n\cdot v_{\mathcal{B}}$.

The arithmeticity of $M_{\mathcal{T}}$ follows from the fact that the fundamental group of $M_{\mathcal{T}}$ is commensurable with the hyperbolic Coxeter group generated by reflections in the facets of the rectified $5$-cell, and the latter group is arithmetic. This follows from the fact that the rectified $5$-cell can be obtained by assembling $320$ copies of the hyperbolic Coxeter simplex $\bar{S}_4 $, which is itself arithmetic, c.f. \cite{JKRT}[p. 342].
\end{proof}

\begin{rem}
If the triangulation $\mathcal{T}$ is orientable, we can lift every isometry between the boundary components of $\mathcal{B}_i$'s to an orientation-reversing isometry of the boundary components of $\widetilde{\mathcal{B}_i}$, and thus obtain the orientable double cover $\widetilde{M_{\mathcal{T}}}$ of the initial manifold $M_{\mathcal{T}}$.
\end{rem}

It follows from the construction that combinatorially equivalent triangulations define isometric manifolds. The converse is also true: given a manifold $M_{\mathcal{T}}$ constructed as described above, we can recover, up to a combinatorial equivalence, the triangulation $\mathcal{T}$. Recall that any manifold of the form $M_{\mathcal{T}}$ is tessellated by a number of copies of the ideal hyperbolic rectified $5$-cell $\mathcal{R}$. 

\begin{prop}\label{prop:recover-triangulation}
The triangulation $\mathcal{T}$  can be uniquely recovered from the decomposition of $M_{\mathcal{T}}$ into copies of $\mathcal{R}$.
\end{prop}

\begin{proof}
Let us introduce an equivalence relation on the copies of $\mathcal{R}$ in the decomposition of $M_\mathcal{T}$, by declaring that $\mathcal{R}_i$ is equivalent to $\mathcal{R}_j$ if they are adjacent along an \textit{octahedral} facet. The equivalence classes naturally correspond to the copies of the block $\mathcal{B}$ which tessellate $M_{\mathcal{T}}$. 

The way in which two copies of the block $\mathcal{B}_k$ and $\mathcal{B}_l$ glue together along their boundary components, is determined by how two adjacent copies of $\mathcal{R}$, say $\mathcal{R}_k\subset \mathcal{B}_k$ and $\mathcal{R}_l\subset \mathcal{B}_l$ are paired along their \textit{tetrahedral} facets. This shows that the decomposition of $M_{\mathcal{T}}$ into copies of the $5$-cell $\mathcal{R}$ allows us to recover the decomposition of $M_{\mathcal{T}}$ into copies of the block $\mathcal{B}$, and this allows us to recover the triangulation $\mathcal{T}$, up to a combinatorial equivalence. 
\end{proof}

Now we want to show that the topology of the hyperbolic manifold $M_{\mathcal{T}}$ solely allows us to recover, up to a combinatorial equivalence, the triangulation $\mathcal{T}$. In order to do so, we need to study the cusp sections of our manifolds.

\subsection{The cusp shape}
Here, we describe the cusp sections of hyperbolic $4$-manifolds defined by triangulations. Let us recall that the cusp shapes of $\mathcal{B}$ are isometric to $K\times I$, where $K$ is the Euclidean Klein bottle depicted in Fig.~\ref{fig:cusps}. When we glue together a number of copies of the block $\mathcal{B}$ to produce the manifold $M_{\mathcal{T}}$, their cusp sections are identified along the boundaries to produce the cusp sections of $M_{\mathcal{T}}$. Each copy of $\mathcal{B}$ contributes its cusp section under these identifications, and together they form a cycle of cusp sections that constitutes a cusp section of $M_{\mathcal{T}}$. Thus, the resulting cusp section of $M_{\mathcal{T}}$ is a closed Euclidean manifold that fibres over $\mathbb{S}^1$, with $K$ as the fibre. The corresponding monodromy is given by an isometry of $K$ into itself, which preserves its tessellation by equilateral triangles shown in Fig.~\ref{fig:cusps}.

\begin{prop}\label{prop:boundary-sequence}
Let $G$ be the group of isometries of $K$ which preserves its tessellation by equilateral triangles. Then $G \cong \mathbb{Z}/2\mathbb{Z}\times \mathfrak{S}_3$.
\end{prop}

\begin{proof}
We begin by showing that there is a short exact sequence 
\begin{equation}
0\rightarrow \mathbb{Z}/2\mathbb{Z}\rightarrow G\rightarrow \mathfrak{S}_3\rightarrow 0.
\end{equation}

Let us notice that the glueing graph $\hat{\hat{\Gamma}}$ for the Klein bottle $K$ is obtained from the graph $\hat{\Gamma}$ by removing three edges which share the same label, as depicted in Fig.~\ref{fig:cuspgraph}.
\begin{figure}[ht]
\centering
\includegraphics[width=0.28\textwidth]{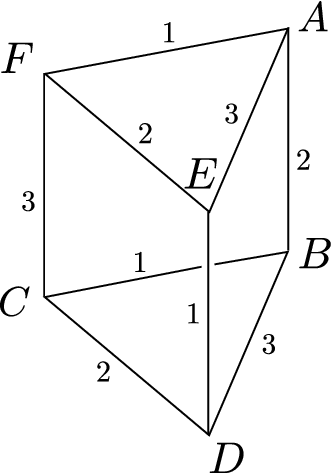}
\caption{The glueing graph $\hat{\hat{\Gamma}}$ for the Klein bottle $K$.}\label{fig:cuspgraph}
\end{figure}

As an unlabelled graph, $\hat{\hat{\Gamma}}$ is isomorphic to the one-skeleton of a triangular prism. There is an order two element in the automorphism group of $\hat{\hat{\Gamma}}$ which exchanges the two bases of the prism, preserving the edge labels. This element induces an involution $P$ of $K$, which is represented by a horizontal reflection of the hexagon in Fig.~\ref{fig:cusps}.

Moreover, any permutation $\sigma\in \mathfrak{S}_3$ of the edge labels is realised by a unique automorphism of $\hat{\hat{\Gamma}}$ as an unlabelled graph, which fixes both triangular bases. Any such automorphism induces an isometry $\psi_{\sigma}$ of $K$. E.g., the transposition of the edges labelled $1$ and $3$ is realised by a reflection in the vertical line of the hexagon in Fig.~\ref{fig:cusps} (the same holds for all other transpositions, up to an appropriate choice of the hexagonal fundamental domain for $K$).

An element of the group $G$ acts on the graph $\hat{\hat{\Gamma}}$ inducing a permutation of the edge labels, and this defines a surjective homomorphism from $G$ onto $\mathfrak{S}_3$. Its kernel is precisely the order two group generated by the involution $P$.

In fact, the above argument shows that the short exact sequence splits. Both factors are normal subgroups of $G$, since one is the kernel of a homomorphism and the second is an index two subgroup. Therefore, the group $G$ decomposes as a direct product: $G \cong \mathbb{Z}/2\mathbb{Z}\times \mathfrak{S}_3$.
\end{proof}

\begin{rem}\label{rem:mapping-class-group}
The group $G$ is naturally mapped into the mapping class group $\mathrm{Mod}\,K$ of the Klein bottle $K$. The isometries $\psi_{\sigma}$, where $\sigma\in \mathfrak{S}_3$ is a transposition, are all isotopic to each other, and therefore they define the same order two element $\tau\in \mathrm{Mod}\,K$. The image of the group $G$ in $\mathrm{Mod}\,K$ is a $\mathbb{Z}/2\mathbb{Z}\times \mathbb{Z}/2\mathbb{Z}$ subgroup generated by $P$ and $\tau$.
\end{rem}

Given a $4$-dimensional triangulation $\mathcal{T}=(\{\Delta_i\}_{i=1}^{2n}, \{g_j\}_{j=1}^{5n})$, let us consider the abstract graph with vertices given by the $20n$ two-dimensional faces of the simplices $\{\Delta_i\}_{i=1}^{2n}$ and edges connecting two vertices if the corresponding two-faces are identified by a pairing map. This graph is a disjoint union of cycles, which we call the \emph{face cycles} corresponding to the triangulation $\mathcal{T}$. 

Associated with every face cycle $c$ of length $h$, there is a sequence 
\begin{equation}
F_0\xrightarrow[\psi_1]\ F_1\xrightarrow[\psi_2]\ \dots \xrightarrow[\psi_{h-1}]\ F_h\xrightarrow[\psi_{h}]\ F_h = F_0
\end{equation}
of triangular faces paired by isometries. This defines a \emph{return map} 
\begin{equation}
r_c=\psi_h\circ\dots\circ\psi_1: F_0 \rightarrow F_0
\end{equation}
as an element of $\mathrm{Symm}\, (F_0) \cong \mathfrak{S}_3$.

From Propositions \ref{prop:isometry-boundary} and \ref{prop:isometry-block} we have the following one-to-one correspondence:
\begin{equation}
\{\mbox{Face cycles of }\mathcal{T}\}\leftrightarrow\{\mbox{Cusps of }M_{\mathcal{T}}\}.
\end{equation}

Moreover, the cusp shape is determined by the length of the associated cycle and its return map, as shown below. 

\begin{prop}\label{prop:cusp-shape}
Let $\mathcal{T}$ be a $4$-dimensional triangulation and let $c$ be a face cycle in $\mathcal{T}$ of length $h$. Let $\mathrm{sign}(h)\in \mathbb{Z}/2\mathbb{Z}$ be the parity of $h$ ($0$ if $h$ is even, and $1$ if $h$ is odd). Let $\mathrm{sign}(r_c)\in \mathbb{Z}/2\mathbb{Z}$ be the parity of the return map $r_c$ as an element of $\mathfrak{S}_3$. The face cycle $c$ defines the element 
\begin{equation}
\phi_c=(P^{\,\mathrm{sign}(h)+\mathrm{sign}(r_c)},r_c)
\end{equation} 
in the group $G \cong \mathbb{Z}/2\mathbb{Z}\times \mathfrak{S}_3$ of automorphisms of the Klein bottle $K$.
The maximal cusp section of the cusp corresponding to the cycle $c$ is isometric to 
\begin{equation}
\frac{K\times [0,h]}{(x,0)\sim (\phi_c(x),h)},
\end{equation}
where $K$ has total Euclidean area of $3\sqrt{3}/2$ (i.e.~the edges of the tessellating equilateral triangles have length one).

\begin{figure}[ht]
\centering
\includegraphics[width=4in]{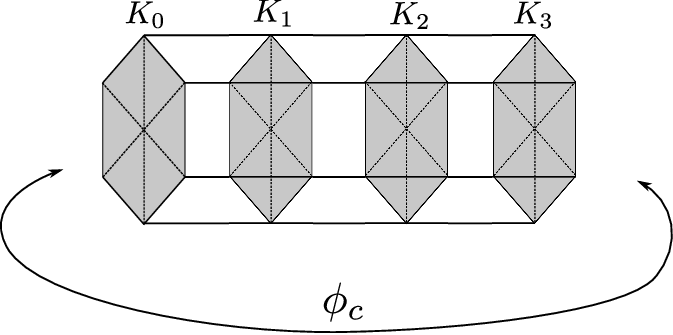}
\caption{The structure of the cusp associated with a cycle of length $3$}\label{fig:cuspmonodromy}
\end{figure}
\end{prop} 

\begin{proof}
Let us observe that the length of the cycle determines the height $h$ of the mapping torus, since we are pairing $h$ blocks of the form $K\times[0,1]$ isometrically along their boundaries. 

Also we notice that, associated with every cycle $c$ of length $h$, there is a sequence of isometries 
\begin{equation}
K_0\xrightarrow[\phi_1]\ K_1\xrightarrow[\phi_2]\ \dots \xrightarrow[\phi_{h-1}]\ K_{h-1}\xrightarrow[\phi_{h}]\ K_h = K_0
\end{equation}
which defines the monodromy 
\begin{equation}
\phi_c=\phi_h\circ\dots\circ\phi_1:K_0\rightarrow K_0.
\end{equation}
 
Clearly, $\phi_c\in G\cong \mathbb{Z}/2\mathbb{Z}\times \mathfrak{S}_3$, since every map $\phi_i: K_i\rightarrow K_{i+1}$ preserves the tessellation by equilateral triangles. The projection of $\phi_c$ onto $\mathfrak{S}_3$ is clearly determined by the return map $r_c$. We need to study the projection of $\phi_c$ onto $\langle P \rangle \cong \mathbb{Z}/2\mathbb{Z}$.

Recall that the group $G$ acts on the glueing graph $\hat{\hat{\Gamma}}$ in Fig.~\ref{fig:cuspgraph}. The element $P$ acts by exchanging its triangular bases, leaving the edge labels unchanged. Thus, the projection of $\phi_c$ onto the group $\langle P \rangle \cong \mathbb{Z}/2\mathbb{Z}$ depends on the behaviour of $\phi_c$ on the bases: it is trivial if it fixes the bases, and non-trivial if it exchanges them.

Recall that there is an inclusion $\hat{\hat{\Gamma}}\subset \hat{\Gamma}$ between the glueing graph of the Klein bottle $K$ and the glueing graph of the boundary manifold $D$ of the block $\mathcal{B}$. The latter is obtained from the former by adding a diagonal to each of the three square faces of the graph $\hat{\hat{\Gamma}}$, as shown in Fig.~\ref{fig:adjacencygraph}, in such a way that the new edges have no common endpoints. There are only two ways of performing this operation. The two resulting labelled graphs correspond to the glueing graphs of two boundary components, say $D_i$ and $D_i'$, of $\mathcal{B}$ which have $K_i$ as one of their cusp sections.

As shown in Fig.~\ref{fig:cuspmonodromy}, each of the maps $\phi_i:K_i\rightarrow K_{i+1}$ can be seen as the composition $p_i\circ a_i$, where
\begin{enumerate}
\item the map $p_i:K_i\times\{1\}\rightarrow K_{i+1}\times\{0\}$ is the restriction to the cusps of a pairing map $g_i$ between the boundary manifolds $D'_i$ and $D_{i+1}$; 
\item the map $a_i$ is an \textit{adjacency} between the boundaries of the cusp sections of the block $\mathcal{B}$, i.e. the map which sends the boundary cusp $K_{i}\times\{0\}$ to $K_{i}\times\{1\}$ according to the rule $(x,0)\rightarrow (x,1)$. 
\end{enumerate}

The adjacency map $a_i: K_{i}\times\{0\}\rightarrow K_{i}\times\{1\}$ can therefore be thought of as acting on the glueing graph $\hat{\Gamma}$ by keeping the glueing graph $\hat{\hat{\Gamma}}$ of $K_{i}$ fixed and exchanging the remaining diagonals, which means taking the glueing graph of the boundary component $D_i$ to the glueing graph of $D'_i$ as depicted in Fig.~\ref{fig:adjacencygraph}.

\begin{figure}[ht]
\centering
\includegraphics[width=0.7\textwidth]{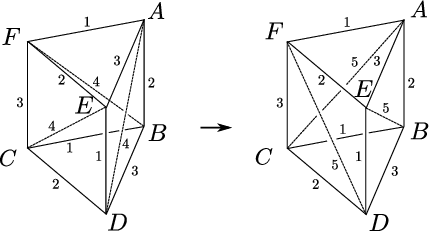}
\caption{Passing from the glueing graph of the boundary component $D_i$ to that of $D'_{i+1}$. The continuous lines represent the initial edges of the glueing graph of $K_i$, while the dashed lines represent the extra edges needed to define boundary components. }\label{fig:adjacencygraph}
\end{figure}

Now we represent each of the maps $\phi_i=p_i\circ a_i:K_i\rightarrow K_{i+1}$ as a self-map of the glueing graph in Fig.~\ref{fig:adjacencygraph} on the left.
Since the map $p_i$ is induced by a map from the glueing graph for $D'_i$ to the glueing graph for $D_{i+1}$, there are two possible cases:
\begin{enumerate}
\item the map $\phi_i$ induces an even permutation of the vertical edges connecting the two triangular bases, and therefore $\phi_i$ exchanges the bases;
\item the map $\phi_i$ induces an odd permutation of the vertical edges, and therefore $\phi_i$ fixes the triangular bases.
\end{enumerate}

Finally, by composing the maps $\phi_i$, $i=1,\dots,h$ we see that
\begin{enumerate}
\item we have an even number of bases exchanges if  either $r_c\in \mathfrak{S}_3$ is odd and the length $h$ of the cycle is odd, or $r_c$ is even and $h$ is even;
\item we have an odd number of bases exchanges if  either $r_c\in \mathfrak{S}_3$ is odd and $h$ is even, or $r_c$ is even and $h$ is odd.
\end{enumerate}
The first case corresponds to the trivial projection of $\phi_c$ onto $\langle P \rangle \cong \mathbb{Z}/2\mathbb{Z}$, while the second case corresponds to a non-trivial one.
\end{proof}

The following fact allows us to recover the cusp shapes of the manifold $M_{\mathcal{T}}$ solely from its topology. 

\begin{prop}\label{prop:cusp-length}
The similarity class of the cusp section associated with a face cycle $c$ determines the maximal cusp section.
\end{prop}

\begin{proof}
In general, the proof follows the main idea of the proof of \cite[Lemma~2.12]{KM}. We provide the reader with a draft of the proof skipping, however, some of the most technical details, which are abundant in this case. 

The cusp sections are endowed with a Euclidean structure, defined up to a similarity transformation. Therefore, any cusp $X$ can be expressed as $X=\mathbb{E}^3/H$, where $H$ is a discrete group of isometries of the Euclidean space $\mathbb{E}^3$.

The maximal cusp section is obtained by choosing for each cusp, which corresponds to a face-cycle of length $h$ in $\mathcal{T}$, the unique section with Euclidean volume equal to $3\sqrt{3}/2\cdot h$. Thus, it suffices to prove that the integer $h$ can be recovered from the geometry of the cusp.

The group $H$ contains a finite-index translation lattice $L<H$, which can be thought of as a lattice in $\mathbb{R}^3$, defined up to a similarity. 

Now let $H^+<H$ be the subgroup of orientation-preserving isometries of $H$. The orientable double cover $\tilde{X}=\mathbb{E}^3/H^+$ of $X=\mathbb{E}^3/H$ fibres over the circle, with the fibre given by the orientable double cover $T$ of $K$. The torus $T$ is represented in Fig.~\ref{fig:orientablecusp}. 
\begin{figure}[ht]
\centering
\includegraphics[width=4in]{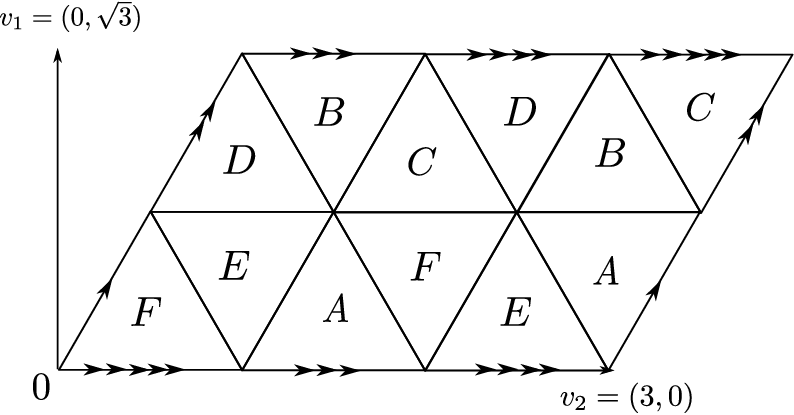}
\caption{The torus $T$, which is the orientable double cover of the Klein bottle $K$, with its Euclidean structure. Up to a homothety, $T$ is generated by translations along the vectors $v_1=(0,\sqrt{3})$ and $v_2=(3,0)$.}\label{fig:orientablecusp}
\end{figure}
Moreover, the monodromy map $\phi:K\rightarrow K$ lifts to a \emph{unique} orientation-preserving isometry of the orientable double cover $T$ of $K$, which we continue to denote by $\phi$.

Notice that we have the inclusions $L<H^+<H$. The length $h$ of the mapping torus remains unchanged while passing to the orientable cover, although the volume of $\tilde{X}$ is twice the volume of $X$.

First, we establish the following auxiliary fact.
 
\begin{lemma}\label{lemma:cusptranslation}
The lifts of elements of $G$ to $T$ are induced by translations of the plane if and only if the return map $r_c\in \mathfrak{S}_3$ is an \emph{even} permutation.
\end{lemma}
\begin{proof}
Notice that, up to a similarity, the torus $T$ from Fig.~\ref{fig:orientablecusp} can be realised as the quotient of the plane $\mathbb{R}^2$ by the translations along the vectors $v_1=(0,\sqrt{3})$ and $v_2=(3,0)$. The orientation-reversing involution $i$ corresponding to the non-trivial automorphism of the cover $T\rightarrow K$ is induced by the isometry 
\begin{equation}
i: (x,y) \mapsto (x+\frac{3}{2},\frac{\sqrt{3}}{2}-y).
\end{equation}

The isometry $P$ exchanges the triangular bases in the glueing graph $\hat{\hat{\Gamma}}$ for the Klein bottle $K$. The isometry of the plane inducing the orientable lift to $T$ is obtained by composing the horizontal reflection 
\begin{equation}
(x,y)\rightarrow (x, \sqrt{3}-y)
\end{equation}
with the isometry inducing the involution $i$. The resulting map is the translation 
\begin{equation}
(x,y)\rightarrow(x+\frac{3}{2},y+\frac{\sqrt{3}}{2}).
\end{equation}
Therefore, the power of $P$ in the monodromy of the cusp, which depends on the length $h$ of the cycle $c$ and the parity of the return map $r_c$,  does not play a role in determining if the monodromy map is a translation.

An odd return map $r_c$ induces an isometry of $K$ corresponding to an odd permutation of the vertices in the bases of $\hat{\hat{\Gamma}}$. The isometry induced by each such $r_c$ lifts to a composition of a reflection along the vertical axis with the isometry of the plane inducing the involution $i$. This is necessarily a point reflection, and therefore not a translation. Moreover, this shows that $r^2_c$ for an odd return map $r_c$, which correspond to an even permutation of the vertices in the bases of $\hat{\hat{\Gamma}}$, is induced by a translation.
\end{proof}

Now we can summarise the dependence of the monodromy map $\phi$ on the length $h$ of the face cycle $c$ and the parity of the return map $r_c$ as follows.

\begin{table}[h]
\begin{center}
\begin{tabular}{|l|p{2.5cm}|p{2.5cm}|}
\hline
& $h$ even & $h$ odd \\ \hline
$r_c$ even & $\phi = (P^0,\tau^0)$ &$\phi = (P^1,\tau^0)$ \\ \hline
$r_c$  odd & $\phi = (P^1,\tau^1)$ &$\phi = (P^0,\tau^1)$ \\
\hline
\end{tabular}\label{tab:monodromy}
\end{center}
\caption{The monodromy map $\phi$ depending on the length $h$ of the face cycle $c$ and its return map $r_c$.}
\end{table}

Now, let $v_1,v_2 \in L$ be two vectors satisfying the following conditions:
\begin{enumerate}
\item[1.] $v_1$ and $v_2$ are orthogonal and the lengths satisfy $l(v_2)=\sqrt{3}\cdot l(v_1)$;
\item[2.] $v_1$ and $v_2$ are the shortest such vectors;
\item[3.] the number $h=\frac{\text{Vol}(\tilde{X})}{l(v_1)^3}$ is an integer.
\end{enumerate}

Notice that the integer $h$ defined above depends only on the similarity class of $X$ and is invariant under rescaling. If such a couple of vectors exists but is not unique, the lengths of the vectors will nonetheless be the same, and therefore the integer $h$ is well defined. Clearly, in the case of the maximal cusp section, the vectors $v_1=(0,\sqrt{3},0)$ and $v_2=(3,0,0)$ satisfy Conditions 1-3. Moreover the integer $h$ is precisely the length of the associated face cycle. We need to show that there never exists a  couple of \emph{shorter} vectors satisfying Conditions 1-3.

Recall that, up to a similarity, we have 
\begin{equation}
\tilde{X}=\frac{T\times [0,h]}{(x,0)\sim (\phi(x),h)}.
\end{equation}
If $h\geq 3$, the vectors $v_1=(0,\sqrt{3},0)$ and $v_2=(3,0,0)$ are clearly the shortest ones satisfying Conditions 1-3.

Now, let us suppose that $h=2$. If the return map $r_c$ of the associated face cycle $c$ is an odd permutation then, by Lemma~\ref{lemma:cusptranslation}, the monodromy map $\phi$ is not induced by a translation, and $L$ is necessarily a proper subgroup of $H^+$, generated by $v_1=(0,\sqrt{3},0)$, $v_2=(3,0,0)$ and a third vector of the form $(n,0,4)$, with $n\in 0,1,2$. In this case we conclude the same as for $h=4$. 

If the return map $r_c$ of the associated face cycle $c$ is an even permutation, the monodromy $\phi$ is induced by a translation. The lattice $L$ is generated by $v_1=(0,\sqrt{3},0)$, $v_2=(3,0,0)$ and a vector of the form $(n,0,2)$, with $n= 0,1,2$. In this case it is clear that the vectors $\pm v_1$ are the shortest non-trivial vectors of the lattice. Therefore the vectors $v_1$ and $v_2$ satisfy Conditions 1-3 again.

Now, let us suppose $h=1$. If the return map $r_c$ on the face cycle is odd then, by Lemma \ref{lemma:cusptranslation}, the monodromy $\phi$ is not induced by a translation. Again, we have a proper inclusion $L<H^+$, and therefore $L$ is generated by $v_1=(0,\sqrt{3},0)$, $v_2=(3,0,0)$ and a vector of the form $(n,0,2)$, with $n\in 0,1,2$. In this case, we conclude as in the previous step.

Finally, we consider the case $h=1$, with the return map $r_c$ given by an even permutation. In this case, the monodromy is induced by a translation. The lattice $L$ is generated by $v_1=(0,\sqrt{3},0)$, $v_2=(3,0,0)$ and a vector $v_3$ of the form $(1/2+n,\sqrt{3}/2,1)$, with $n=0,1,2$.

In the case $n=0$, a computation with a SAGE routine \cite{K} (see also Appendix B) shows that the only non-trivial vectors $v\in L$ that satisfy $l(v)\leq 3$ and Condition 3 are $\pm v_1$ and $\pm v_2$. 

If $n=1$ or $n=2$, then the same SAGE routine shows that the pair of shortest vectors in $L$ satisfying Conditions 1-3 are again $\pm v_1$ and $\pm v_2$. In fact, both cases produce isometric lattices, up to a reflection in the horizontal axis.
\end{proof}

\begin{prop}\label{prop:Epstein-Penner}
Given a $4$-dimensional triangulation $\calT$, the Epstein-Penner decomposition of the manifold $M_{\calT}$ with respect to the maximal cusp section is the decomposition into copies of the rectified $5$-cell $\calR$.
\end{prop}

\begin{proof}
The maximal cusp section and the tessellation of the manifold $M_{\calT}$ into copies of the rectified $5$-cell $\calR$ lift to a tessellation of $\mathbb{H}^4$ into copies of $\calR$, together with a horocusp at each ideal vertex. The set of horocusps is invariant under the isometry group of the tessellation. The Epstein-Penner decomposition \cite{EP} is obtained by interpreting the horocusps as points on the light cone in $\mathbb{R}^{4,1}$, and taking their Euclidean convex hull. Thus, by symmetry, the resulting decomposition is just the original decomposition into rectified $5$-cells.
\end{proof}

\begin{cor}\label{cor:buildingblock}
Any isometry of the building block $\cal{B}$ has to preserve its decomposition into copies of the rectified $5$-cell $\cal{R}$.
\end{cor}

\begin{proof}
Consider the manifold $\calB^{\prime}$ obtained by doubling the building block $\cal{B}$ in its boundary, and lift the tesselation of $\calB$ into copies of $\calR$ to a tessellation of $\calB^{\prime}$. The manifold $\calB^{\prime}$ is associated to the triangulation obtained by doubling a $4$-dimensional simplex in its boundary. 

Any isometry of the building block $\calB$ induces a unique isometry of its double $\calB^{\prime}$ which fixes each of the two copies of $\calB$ and, by Proposition \ref{prop:Epstein-Penner}, has to preserve its decomposition into copies of $\cal{R}$. Therefore any isometry of $\calB$ preserves its decomposition into copies of $\calR$.
\end{proof}

\begin{teo}
The topology of the manifold $M_{\mathcal{T}}$ is determined by the triangulation $\mathcal{T}$, up to a combinatorial equivalence, and vice versa.
\end{teo}

\begin{proof}
The topology of $M_{\mathcal{T}}$ determines its hyperbolic structure uniquely, grace to the Mostow-Prasad rigidity. The hyperbolic structure, in its own turn, determines the similarity class of the cusp sections. By Proposition \ref{prop:cusp-length}, the cusp shapes determine the maximal cusp section. The latter determines the decomposition of $M_{\mathcal{T}}$ into copies of the rectified $5$-cell $\mathcal{R}$, which is, by Proposition \ref{prop:Epstein-Penner}, the Epstein-Penner decomposition corresponding to the maximal cusp section. Thus, the triangulation $\mathcal{T}$ can be recovered, up to a combinatorial equivalence, according to Proposition \ref{prop:recover-triangulation}.
\end{proof}

\begin{teo}\label{teo:aut-isom}
The group $\mathrm{Aut}\, \mathcal{T}$ of combinatorial equivalences of a $4$-dimensional triangulation 
\begin{equation}
\mathcal{T}=(\{\Delta_i\}_{i=1}^{2n}, \{g_j\}_{j=1}^{5n})
\end{equation}
and the group $\mathrm{Isom}\,M_{\mathcal{T}}$ of isometries of the associated manifold $M_{\mathcal{T}}$ are isomorphic. 
\end{teo}

\begin{proof}
Every simplicial map $\phi_{kl}:\Delta_k\rightarrow \Delta_l$ determines an isometry between the corresponding copies $\mathcal{B}_k$ and $\mathcal{B}_l$ of the building block, according to Proposition \ref{prop:isometry-block}. By applying these isometries to each copy of the building block, we obtain an isometry of $M_{\mathcal{T}}$. This defines a homomorphism from $\mathrm{Aut}\,\mathcal{T}$ to  $\mathrm{Isom}\,M_{\mathcal{T}}$.

To build an inverse of this homomorphism, we notice that every isometry of $M_{\mathcal{T}}$ has to preserve the Epstein-Penner decomposition into copies of the rectified $5$-cell $\mathcal{R}$. Therefore, following the proof of Proposition \ref{prop:recover-triangulation}, we see that every isometry has to preserve the decomposition into copies of the building block $\mathcal{B}$, and thus it defines a combinatorial equivalence of the triangulation $\mathcal{T}$ to itself, which is an element of $\mathrm{Aut}\,\mathcal{T}$.
\end{proof}

Now we can reprove some result initially obtained in \cite{KM} and \cite{S2}:
\begin{teo}
There exists a hyperbolic $4$-manifold with one cusp. Such a manifold can be chosen to be both orientable and non-orientable. Its cusp section is, respectively, either a three-dimensional torus or a product $S^1\times K$, where $K$ is a Klein bottle.
\end{teo}

\begin{proof}
Let $A$ and $B$ be two $4$-dimensional simplices. Let $1$, $2$, $\dots$, $5$ denote the vertices of $A$, and let $\bar{1}$, $\bar{2}$, $\dots$, $\bar{5}$ denote those of $B$. The following facet pairing defines a triangulation with a single face cycle $c$, which has even length $h=20$ and whose associated return map $r_c$ is an even permutation. 
\begin{eqnarray}
(1234)\mapsto (\overline{4231}),\\
(1235)\mapsto (\overline{3521}),\\
(2345)\mapsto (\overline{3254}),\\
(1345)\mapsto (\overline{4513}),\\
(1245)\mapsto (\overline{5421}).
\end{eqnarray}
E.g., if we start from the triangular face $(123)$, then $r_c$ is the identity map $(123) \mapsto (123)$. Thus, the triangulation $\mathcal{T}$ that consists of the simplices $A$ and $B$ with the above facet pairings produces a non-orientable hyperbolic $4$-manifold $M := M_{\mathcal{T}}$ with one cusp. The cusp shape is homeomorphic to $S^1\times K$. By taking the orientable double-cover $\widetilde{M}$ of $M$, we obtain an orientable hyperbolic $4$-manifold with one cusp. Indeed, the cusp section of $M$ is non-orientable, and thus lifts to a single connected component of the cusp section of $\widetilde{M}$. The cusp shape of $\widetilde{M}$ is a three-dimensional torus, according to Table~\ref{tab:monodromy}. The volumes of $M$ and $\widetilde{M}$ are respectively $2\cdot v_{\mathcal{B}} = 8\,\pi^2/3$ and $4\cdot v_{\mathcal{B}} = 16\,\pi^2/3$.
\end{proof}

\begin{quest}
Does there exist a hyperbolic $4$-manifold with a single cusp and volume $< 8\,\pi^2/3$ ($< 16\,\pi^2/3$ in the orientable case)?
\end{quest}

\section{Symmetries of a triangulation}

Below, for any given finite group $G$, we describe a construction of an orientable triangulation $\mathcal{T}$ such that $\mathrm{Aut}\,\mathcal{T} \cong G$. Then, by applying Theorem~\ref{teo:aut-isom}, we can produce two complete hyperbolic $4$-manifolds of finite volume: a non-orientable manifold $\mathcal{M}$, such that $\mathrm{Isom}\,\mathcal{M}\cong G$, and an orientable manifold $\widetilde{\mathcal{M}}$ such that $\mathrm{Isom}^+\,\widetilde{\mathcal{M}}\cong G$. Finally, we estimate the volume of our manifolds in terms of the order of the group $G$.

\subsection{Orientation of a simplex}

Let $S_4$ be the regular $4$-simplex with vertices labelled $v_1$, $v_2$, $\dots$, $v_5$. Also, each facet of $S_4$ gets a label: the label of its opposite vertex with respect to the natural self-duality of $S_4$. The orientation on $S_4$ is defined by some order $v_{i_1} < v_{i_2} < \dots < v_{i_5}$ on its vertices, which we also denote by $[v_{i_1}, v_{i_2}, \dots, v_{i_5}]$.

Let, up to a suitable change of notation, $(v_1, v_2, \dots, v_5)$ be the standard oriented $4$-simplex and $[v_1, v_2, \dots, v_5]$ be its positive orientation. The orientation of its $3$-dimensional facet is obtained from the classical formula for the boundary of a simplex:
\begin{equation}\label{eq:orientation}
\partial [v_1, v_2, \dots, v_5] = \sum^5_{i=1} (-1)^{i+1}\, [v_1, \dots, \widehat{v_i}, \dots, v_5],
\end{equation} 
where the hat sign means that the respective vertex is omitted. Thus, the orientation of a facet of $S_4$ either coincides with that induced by the order on its vertices, or is opposite to it.

Let now $A$ and $B$ be two copies of $S_4$, and $\phi_{v_i u_j}$ be a map identifying a pair of their facets labelled $v_i$ and $u_j$, correspondingly. Then, given the orientations of $A$ and $B$, and subsequently those of each facet from formula \eqref{eq:orientation}, we can easily determine if $\phi_{v_i u_j}$ is orientation preserving or reversing.

\subsection{Constructing a triangulation}

Let $\mathcal{T}$ be a $4$-dimensional triangulation. We denote by $\|\mathcal{T}\|$ the number of simplices in $\mathcal{T}$. Let $\mathrm{Aut}\,\mathcal{T}$ denote the automorphism group of $\mathcal{T}$. If $\mathcal{T}$ is orientable, let $\mathrm{Aut}^{+}\,\mathcal{T}$ denote the subgroup of orientation-preserving automorphisms of $\mathcal{T}$. 

\begin{teo}\label{teo:aut-triangulation}
For any finite group $G$ of rank $m$ with $n$ elements, there exists an orientable triangulation $\mathcal{T}$, such that $\mathrm{Aut}^{+}\,\mathcal{T} \cong \mathrm{Aut}\,\mathcal{T} \cong G$. Also, $\|\mathcal{T}\| \leq C\cdot n\cdot m^2$, for some constant $C>1$ which does not depend on $G$. 
\end{teo}

Before starting the proof of Theorem~\ref{teo:aut-triangulation}, we need to produce some more technical ingredients. First, we produce a ``partial'' $4$-dimensional triangulation which consists of two $4$-simplices $A$ and $B$, has two unpaired facets and admits no non-trivial automorphisms. We call it the \textit{edge sub-complex}, since later on it will be associated with some edges in the Cayley graph of the group $G$. 

\begin{figure}[ht]
\centering
\includegraphics[width=2in]{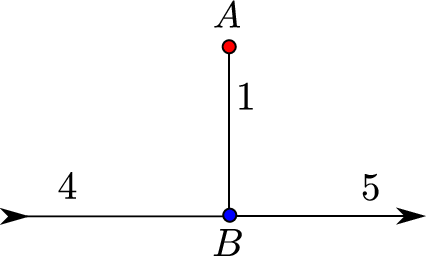}
\caption{A schematic picture of the edge sub-complex $E$.}\label{fig:complexE}
\end{figure}

Let us proceed as follows: take two $4$-simplices $A$ and $B$, label the vertices of both with the numbers $1$, $2$, $\dots$, $5$. Their facet also become labelled, by the self-duality of $S_4$. The orientation of $A$ is $[1, 2, 3, 4, 5]$, and that of $B$ is the opposite, say $[5, 2, 3, 4, 1]$. First, we identify four facet of $A$ in pairs. Namely,
\begin{eqnarray}\label{eq:vertex-block-ident-a}
\mbox{facets $2$ and $3$ by the map } (1345) \longleftrightarrow (4251);\\
\mbox{facets $5$ and $4$ by the map } (1234) \longleftrightarrow (2351).
\end{eqnarray}
Let us notice that both maps above are orientation reversing. 

The simplicial complex defined by pairing the facets of $A$ according to the maps \eqref{eq:vertex-block-ident-a} has trivial automorphism group. Indeed, any of its automorphisms $\phi$ has to be induced by an automorphism $\tilde{\phi}: S_4 \longrightarrow S_4$ of the simplex $A$. Since the facet labelled $1$ is unpaired, it is taken by $\phi$ into itself. Thus, $\tilde{\phi}$ fixes the vertex labelled $1$. Moreover, any automorphism $\phi$ preserves the facet pairing. Thus, in our case, $\tilde{\phi}$ has to preserve the sets $\{2, 3\}$ and $\{4, 5\}$. Since $\tilde{\phi}\in \mathfrak{S}_5$, it can be either the identity, or any of the transpositions $(23)$ and $(45)$, or their product. By a straightforward calculation, we verify that, apart from the identity, in each other case at least one facet pairing in \eqref{eq:vertex-block-ident-a} is not preserved. 

Now we pair the facet labelled $1$ of $A$ with the facet labelled $1$ of $B$ by the identity map. However, since the orientations of $A$ and $B$ are opposite, this map is orientation-reversing, as well. Finally, we identify the following facets of $B$:  
\begin{equation}\label{eq:vertex-block-ident-b}
\mbox{facets $2$ and $3$ by the map } (1345) \longleftrightarrow (4251).
\end{equation}
This map is also orientation reversing. Thus, we obtain an orientable triangulation with two unidentified facets.

We show that the resulting simplicial complex has trivial automorphism group. In this case, the simplex $A$ has no unpaired facets, and $B$ has two unpaired ones. Thus, no automorphism $\phi$ of the resulting complex can exchange its simplices $A$ and $B$. In this case, $\phi$ acts trivially on $A$. Since $\phi$ preserves the facet identification between facet $1$ of $A$ and facet $1$ of $B$, which is the identity map, $\phi$ acts also trivially on $B$.

Let us call the resulting complex $E$. It is schematically depicted in Fig.~\ref{fig:complexE}. The simplices $A$ and $B$ are represented by vertices, the identifications of the same simplex facet are not shown, the edge between $A$ and $B$ represents the pairing of facet $1$ of $A$ with facet $1$ of $B$. The unidentified facets of $B$ are represented by the edges labelled respectively $4$ and $5$. The arrows emphasise the fact, that there is no automorphism exchanging the unpaired facets of $E$. 

\begin{figure}[ht]
\centering
\includegraphics[width=0.8\textwidth]{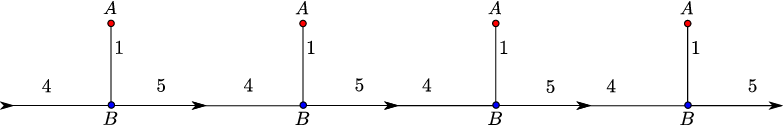}
\caption{A schematic picture of the edge complex $E_n$.}\label{fig:complexEn}
\end{figure}

We link consequently $n\geq 1$ copies of $E$ in order to obtain a simplicial complex $E_n$, with $E_1 = E$. This construction is schematically depicted in Fig.~\ref{fig:complexEn}. The identification between the respective facets labelled $4$ and $5$ of two consequent copies of $E$ is given by
\begin{equation}
(1235) \longleftrightarrow (1234).
\end{equation}
This facet pairing is orientation reversing. Thus, the complex $E_n$ is an orientable triangulation with two unpaired facets. We call $E_n$, $n\geq 1$, the \textit{edge complexes}, although we have a family of mutually non-isomorphic complexes. Analogous to the above, we can show that each $E_n$ has trivial automorphism group.

\medskip

\prf{Theorem~\ref{teo:aut-triangulation}} Given a finite group $G$, let us describe an orientable triangulation $\mathcal{T}$, such that $\mathrm{Aut}\,\mathcal{T}$ is isomorphic to $G$. Let $G$ be of order $n$ and let it have a presentation with $m$ generators $\{s_1, s_2, \dots, s_m\}$ and a certain number of relations. Let $\Gamma$ be the Cayley graph of $G$ corresponding to that presentation. Each edge of $\Gamma$ is directed: it connects a vertex labelled by $g\in G$ with one labelled by $g\cdot s_i$. Here we stress the fact that in the case of an order two generator $s_i$, we introduce \emph{two} oriented edges between every couple of vertices $\{g, g\cdot s_i\}$: the first has $g$ as starting point and $g\cdot s_i$ as endpoint, while the other is oriented in the opposite way.
Also, each edge $[g, g\cdot s_i]$ is labelled by the corresponding generator $s_i$. Considering $\Gamma$ as a directed labelled graph, we have exactly $\mathrm{Aut}\,\Gamma \cong G$, see e.g. \cite{Meier}. 

\begin{figure}[ht]
\centering
\includegraphics[width=2in]{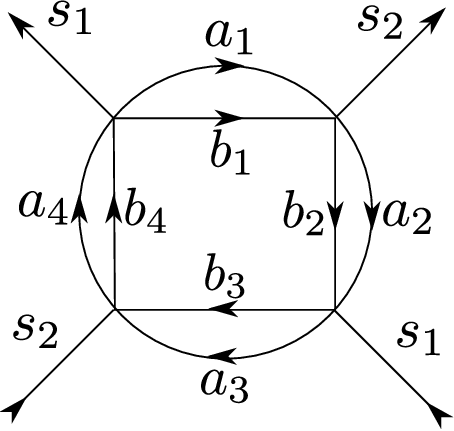}
\caption{A vertex of the graph $\Gamma$ blown up, in the case $m=2$.}\label{fig:vertexBlowUp}
\end{figure}

Given the graph $\Gamma$, we ``blow-up'' each vertex as depicted in Fig.~\ref{fig:vertexBlowUp}, and add additional labels $a_i$, $b_i$ on the new edges. Let us call the modified graph $\Gamma^{\prime}$. Now the graph $\Gamma^{\prime}$ is uniformly $5$-valent, and still $\mathrm{Aut}\,\Gamma^{\prime}\cong G$, as a group of automorphisms of a directed labelled graph. More formally, we perform the following on the graph $\Gamma$ in order to produce $\Gamma^{\prime}$:
\begin{itemize}
\item[1.] fix an order on the generating set of $G$: $s_1 < s_2 < \dots < s_m$;
\item[2.] replace every vertex $g\in \Gamma$ with $2m$ vertices labelled
\begin{equation}
(g, s_1, -), (g, s_1, +), \dots,(g, s_m, -), (g, s_m, +);
\end{equation}
\item[3.] for each edge $[g, g\cdot s_i]$ in the initial graph $\Gamma$, we add an edge connecting $(g, s_i, -)$ to $(g\cdot s_i, s_i, +)$;
\item[4.] for each vertex $g\in G$ we add $4m$ directed edges labelled by some new labels $a_i$, $b_i$, $i=1, \dots, 2m$, such that each edge labelled $a_i$ or $b_i$, $i=1, \dots, m-1$ connects $(g, s_i, -)$ to $(g, s_{i+1}, -)$, each edge labelled $a_i$ or $b_i$, $i=m+1, \dots, 2m$, connects $(g, s_i, +)$ to $(g, s_{i+1}, +)$, the edge labelled $a_m$ or $b_m$ connects $(g, s_m, -)$ to $(g, s_1, +)$, the edge labelled $a_{2m}$ or $b_{2m}$ connects $(g, s_m, +)$ to $(g, s_1, -)$; 
\item[5.] finally, we extend the order on the generating set of $G$ to an order on all edge labels:
\begin{equation}
s_1 < s_2 < \dots < s_m < a_1 < \dots < a_{2m} < b_1 < \dots < b_{2m}.
\end{equation}
\end{itemize}
We obtain a directed labelled graph $\Gamma^{\prime}$, such that $\mathrm{Aut}\,\Gamma^{\prime} \cong \mathrm{Aut}\,\Gamma \cong G$. To see that the automorphism group of $\Gamma^{\prime}$ is the same as the automorphism group of $\Gamma$, notice that the subgraphs which result from the blow-up of the vertices of $\Gamma$ can be recovered by the combinatorial properties of the graph $\Gamma^{\prime}$ as follows.

Let us call two vertices $v$ and $w$ of the graph $\Gamma^{\prime}$ \textit{related} if there are two edges connecting them, and extend this to an equivalence relation. The equivalence classes are in a one-to-one correspondence with the sub-graphs that result from the blow-up of the vertices of $\Gamma$ (\emph{i.e.} with the elements of the group $G$). Since the number of edges connecting two vertices is preserved under any automorphism of $\Gamma^{\prime}$, the equivalence classes are necessarily permuted under the action of $\text{Aut}\,\Gamma^{\prime}$. As a consequence, any directed, labelled automorphism of $\Gamma^{\prime}$ induces an automorphism of $\Gamma$. The opposite correspondence is obvious.

The graph $\Gamma^{\prime}$ is uniformly $5$-valent and has a total order on its edge labels. Now we associate a triangulation $\mathcal{T}$ with it.
First, with each vertex $v\in \Gamma^{\prime}$ we associate a simplex $S_v$ such that its vertices inherit respectively the labels of the edges adjacent to $v$. Thus, the vertices of $S_v$ have an order induced by that on the edge labels of $\Gamma^{\prime}$. Namely,
\begin{itemize}
\item[1.] the simplex associated with the vertex $(g, s_i, -)$, $i\neq 1$, has vertices
\begin{equation}
s_i < a_{i-1} < a_i < b_{i-1} < b_i;
\end{equation}
\item[2.] the simplex associated with the vertex $(g, s_i, +)$ has vertices
\begin{equation}
s_i < a_{m+i-1} < a_{m+i} < b_{m+i-1} < b_{m+i};
\end{equation}
\item[3.] the simplex associated with the vertex $(g, s_1, -)$ has vertices
\begin{equation}
s_1 < a_1 < a_{2m} < b_1 < b_{2m}.
\end{equation}
\end{itemize}

We assume that the orientation of each simplex coincides with the orientation induced by its vertex labels. Again, for each simplex $S_v$, its facets are labelled with the same labels as its vertices in accordance with self-duality: a facet has the label of the opposite vertex.

Second, we produce $5m$ non-isomorphic copies of edges complexes $E_{s_i}$, $E_{a_i}$ and $E_{b_i}$, such that each complex is indexed by the respective edge label of $\Gamma^{\prime}$: we put $E_{s_i} := E_i$, $E_{a_i} := E_{m+i}$, and $E_{b_i} := E_{3m+i}$. Notice that we can make this choice in such a way that each of the chosen copies of the edge sub-complex has less than $10\cdot m$ simplices. Then, for each directed edge of $\Gamma^{\prime}$ labelled $a$ and joining two vertices $u$ and $v$, we glue in the edge sub-complex $E_a$, in such a way that 
\begin{itemize}
\item[1.] the direction of the arrows in $E_a$ coincides with the direction of the edge $a$;
\item[2.] the facets of $E_a$ labelled $4$ and $4$, previously unidentified, are paired respectively with the facets of $S_u$ and $S_v$ labelled $a$. 
\end{itemize}

Above, if we identify two facets (one of $E_a$, another of $S_u$ or $S_v$) with distinct orientations, we define the pairing map in such a way that it preserves the order of the vertices. If we identify two facets whose orientations are the same (both either coincide with the orientation given by the vertex order, or both are opposite to it), we choose a pairing map preserving the order of the vertices, and compose it with the restriction of the orientation-reversing map $(12345) \longleftrightarrow (13245)$ to the appropriate facet of $E_a$ (which is necessarily labelled $4$ or $5$ by construction). Thus, our triangulation $\mathcal{T}$ can be made orientable. 

Any automorphism $\phi$ of $\mathcal{T}$ has to preserve the number of self-pairings of each simplex in it, hence each edge sub-complex $E_a$ is mapped onto its own copy contained in $\mathcal{T}$. Also, the arc directions in the edge sub-complexes and their labels are preserved by $\phi$. Thus, $\phi$ induces an automorphism $\tilde{\phi}$ of the directed labelled graph $\Gamma^{\prime}$. Since the edge complexes have trivial automorphism group, an element $\tilde{\phi} \in \mathrm{Aut}\,\Gamma^{\prime}$ induces a \emph{unique} automorphism of the triangulation $\mathcal{T}$. Thus $\mathrm{Aut}\,\mathcal{T} \cong \mathrm{Aut}\,\Gamma^{\prime} \cong \mathrm{Aut}\, \Gamma$.

Moreover, suppose that an automorphism $\phi\in \mathrm{Aut}\,\mathcal{T}$ is orientation-reversing. Then it reverses the orientation of each simplex in some edge sub-complex $E_a$. By the construction of $E_a$, such a complex is \textit{chiral}, meaning that inverting the orientation of each simplex in $E_a$ produces a complex $\overline{E_a}$, which is not isomorphic to $E_a$: the direction of the arrows in $E_a$ becomes inverse. The contradiction shows that indeed $\mathrm{Aut}^{+}\,\mathcal{T} \cong \mathrm{Aut}\, \mathcal{T}$. 

Finally, we compute the number of simplices $\|\mathcal{T}\|$ in the triangulation $\mathcal{T}$. Given that the order of the group $G$ is $n$, and that we choose a generating set with $m$ elements, we have the following amount of simplices in $\mathcal{T}$:
\begin{itemize}
\item[1.] one simplex for each vertex of the modified Cayley graph $\Gamma^{\prime}$ of $G$,
\item[2.] some amount of simplices in each edge sub-complex for each edge of $\Gamma^{\prime}$.
\end{itemize}

We observe that by construction the former amount above grows like $C_1\cdot n\cdot m$, for some $C_1 > 1$ independent of the group $G$, since we have $2\cdot n \cdot m$ vertices in the graph $\Gamma^{\prime}$ after ``blowing-up'' the initial Cayley graph $\Gamma$ of $G$ having $n$ vertices. The latter amount of simplices above grows like $C_2\cdot n\cdot m^2$, with some $C_2 > 1$ independent of $G$, since we have $5\cdot m\cdot n$ edges in $\Gamma^{\prime}$, and each edge corresponds to an edge sub-complex with $\leq 10\cdot m$ simplices in each. Finally, we obtain the desired estimate for $\|\mathcal{T}\|$. \prfend

By an observation of Frucht \cite{Frucht}, a finite group with $n$ elements has rank at most $m \leq \log n/\log 2$.  Combining this with Theorem~\ref{teo:aut-isom} and Theorem~\ref{teo:aut-triangulation}, we obtain the following corollary.

\begin{cor}
Given a finite group $G$ with $n$ elements there is a hyperbolic $4$-manifold $\mathcal{M}$, such that $\mathrm{Isom}\,\mathcal{M} \cong G$ and $\mathrm{Vol}\,\mathcal{M} \leq C\cdot n\cdot \log^2 n$, for some $C>1$ independent of $G$. 
\end{cor}

By passing to the orientable double-cover of the above manifold, we obtain one more corollary below.

\begin{cor}
Given a finite group $G$ with $n$ elements there is an orientable hyperbolic $4$-manifold $\mathcal{M}$, such that $\mathrm{Isom}^+\,\mathcal{M} \cong G$ and $\mathrm{Vol}\,\mathcal{M} \leq C\cdot n\cdot \log^2 n$, for some $C>1$ independent of $G$.
\end{cor}

Given a finite group $G$ and $n\geq 2$, define as in \cite{BL} $f(n,G)$ to be the minimal volume of a hyperbolic $n$-manifold $\mathcal{M}$ such that $ \mathrm{Isom}\,\mathcal{M} \cong G$. We have the following corollary:

\begin{cor}\label{cor:volume-growth}
There exist constants $C_1$ and $C_2$ such that $$ C_1\cdot |G| < f(4,G) < C_2\cdot |G| \cdot \log^2(|G|).$$
\end{cor}

\begin{proof}
The first inequality is an easy consequence of the Kazhdan-Margulis theorem \cite{KMa}. The second inequality comes from the construction above.
\end{proof}

Concerning the above corollaries, we should mention that there are examples of hyperbolic $4$-manifolds with sufficiently large isometry group and relatively small volume \cite[Table~5.1]{FT2014}, which cannot be produced by the construction of Theorem~\ref{teo:aut-triangulation}.

Finally, it is worth mentioning an additional property of the manifolds that we have constructed:
\begin{prop}
Let $G$ be a finite group and let $\Gamma$ be its Cayley graph with respect to a given presentation. Let $\mathcal{M}_{\Gamma}$ be a manifold such that $\mathrm{Isom}\, \mathcal{M}_{\Gamma} \cong G$, constructed as in the proof of Theorem \ref{teo:aut-triangulation}. Then the isometric action of $G$ on 
$\mathcal{M}_{\Gamma}$ is free, and the quotient $\calM_{\Gamma}/G$ is a hyperbolic $4$-manifold with trivial isometry group.
\end{prop}

\begin{proof}
By our construction, is sufficient to check that a non-trivial element $G$ acts on the triangulation associated with $\mathcal{M}_{\Gamma}$ without fixing any simplex and without exchanging any two different simplices which are paired together along a facet. 

The first property is a consequence of the fact that a non-trivial element of $G$  acts on its Cayley graph $\Gamma$ without fixing any vertex or edge. 

The number of self pairings between facets of a simplex is clearly preserved under the action of the group $G$. By construction, whenever two different simplices of the triangulation are paired together along a facet they have a different number of self-pairings, therefore they cannot be exchanged under the action of a non-trivial element of $G$.
\end{proof}

\begin{rem}
The construction introduced in the proof of Theorem \ref{teo:aut-triangulation} depends only on the \emph{local} properties of the Cayley graph of the group $G$. As a consequence, it is easily adaptable to the case of \emph{any} finitely generated group, proving that such a group is the isometry group of a complete hyperbolic $4$-manifold, possibly of infinite volume.

\end{rem}

\subsection{Manifolds with given symmetries}

In this section we describe a construction that produces a family of pairwise non-isometric complete finite-volume hyperbolic four-manifolds $M$ with $\mathrm{Isom}\,M \cong G$ and volume $\mathrm{Vol}\,M \leq V$ for a given finite group $G$ and some $V \geq V_0 > 0$, big enough. The number of manifolds $\rho_G(V) = \#\{\, M\, |\, \mathrm{Isom}\,M \cong G \mbox{ and } \mathrm{Vol}\,M \leq V \}$ in this family grows super-exponentially with respect to $V$. More precisely, we shall prove the following statement:

\begin{teo}\label{teo:superexponential}
For any finite group $G$, there exists a $V_0 > 0$ sufficiently large such that for all $V \geq V_0$ we have
$\rho_G(V) \geq C^{\, V \log V}$, for some $C>1$ independent of $G$.
\end{teo}

The idea of the proof is as follows: given a finite group $G$, we build many non-equivalent edge complexes, ``modelled'' on trivalent graphs. This allows us to build a class of combinatorially non-equivalent triangulations with $G$ as automorphism group. In order to estimate the number of such triangulations we use the following result by Bollob\'{a}s \cite{Bo}:

\begin{prop}\label{prop:trivalent-number}
Let $g(k)$ denote the number of unlabelled trivalent graphs on $2 k$ vertices. There exists a constant $C_0>1$ such that the function $g(k)$ behaves asymptotically like $C_0^{\,k\log{k}}$ for $k\to\infty$: \begin{equation}
g(k)\sim C_0^{\,k\log{k}}.
\end{equation}
\end{prop}

Since we want our complexes to have trivial automorphism group, it is convenient to restrict our attention to \textit{asymmetric} graphs, i.e. graphs with trivial automorphism group. It turns out \cite{Bo, KSV} that there are plenty of asymmetric graphs: 

\begin{prop}\label{prop:trivalent_asymmetric}
Let $g(k)$ denote the number of unlabelled trivalent graphs on $2 k$ vertices, and let $f(k)$ denote the number of asymmetric graphs among them. Then 
\begin{equation}\lim_{k\to \infty}\frac{f(k)}{g(k)} = 1. 
\end{equation}
\end{prop}

By combining the above two results, it is easy to see that the number $f(k)$ of asymmetric trivalent graphs on $2 k$ vertices grows like $C_0^{\,k\log{k}}$:
\begin{equation}\label{eq:asymptotic-graphs}
f(k)\sim C_0^{\,k\log{k}}.
\end{equation}

\prf{Theorem~\ref{teo:superexponential}} Let $\Gamma$ be the Cayley graph of $G$, a finite group of order $n$ and rank $m$. Let $\Gamma^\prime$ be the modified Cayley graph, as in the proof of Theorem \ref{teo:aut-triangulation}. The number of vertices of the graph $\Gamma^\prime$ equals $v = 2\cdot n\cdot m$, and the number of edges $e$ satisfies the relation $2\cdot e = 5\cdot v$. Thus, $e = 5\cdot n\cdot m$. Moreover, recall that the edges of $\Gamma^\prime$ are labelled with $5m$ edge labels, and the group of labelled automorphisms of $\Gamma^\prime$ is isomorphic to $G$.

Analogous to the proof of Theorem \ref{teo:aut-triangulation}, we need to build $5 m$ orientable, pairwise non-isomorphic four-dimensional simplicial complexes with two unpaired facets and trivial automorphism group.
In order to do so, let us choose an integer $N$ such that $f(N)\geq 5 m$, and let us choose $5m$ asymmetric trivalent graphs $\Gamma_1,\dots,\Gamma_{5m}$ on $2 N$ vertices.

Let $\Gamma_i$ be such a graph. We modify it by deleting one of its edges (such that it remains connected, e.g. any edge in the complement of a spanning tree), and leaving two vertices of valence two. Let $\Gamma^\prime_i$ be the resulting graph. If $\Gamma_i$ and $\Gamma_j$ are not isomorphic, then $\Gamma^\prime_i$ and $\Gamma^\prime_j$ are not isomorphic either. Suppose the contrary: there exists an isomorphism $\phi: \Gamma^\prime_i \rightarrow \Gamma^\prime_j$. Then $\phi$ takes the only pair of valence two vertices of $\Gamma^\prime_i$ to the respective pair of valence two vertices in $\Gamma^\prime_j$. By bringing back the edge connecting these vertices, we obtain an isomorphism $\widetilde{\phi}: \Gamma_i \rightarrow \Gamma_j$, which is a clear contradiction. Moreover the graphs $\Gamma^\prime_i$ are still asymmetric, since any automorphism $\psi$ of $\Gamma^\prime_i$ extends to an automophism  $\widetilde{\psi}$ of $\Gamma_i$. Since  $\Gamma_i$ is asymmetric, $\widetilde{\psi}$ is necessarily the identity, therefore $\psi$ is the identity.

Now, let $S_v$ denote an oriented four-simplex for each vertex $v$ of $\Gamma^\prime_i$, two of whose facets are identified by an orientation-reversing map. Suppose that the vertices of $S_v$ are labelled $1$ to $5$ and the orientation of $S_v$ is $[1,2,3,4,5]$. We perform the identification $(1234) \longleftrightarrow (5342)$, while all other facets remain unpaired. Let us call $S_v$ a \emph{vertex sub-complex}.

Now, let $e_1,\dots,e_{3 N - 1}$ be the edges of $\Gamma^\prime_i$. With each edge $e_i$, we associate a copy $E^\prime_i$ of the edge sub-complex $E$ defined in the proof of Theorem \ref{teo:aut-triangulation} and represented in Fig.\ \ref{fig:complexE}. Whenever an edge $e_i$ of $\Gamma^\prime_i$ connects two vertices $v$ and $w$, we identify with an orientation reversing map a yet unpaired facet of the vertex sub-complex $S_v$ with an unpaired facet of the edge sub-complex $E^\prime_i$, and a yet unpaired facet of the vertex sub-complex $S_w$ with the remaining unpaired facet of $E^\prime_i$. Notice that there are many possible choices, both in which facet of a vertex sub-complex we choose to pair and in which pairing maps we choose. We denote the resulting complex, which has two unpaired facets, by $C_i$.

The $5m$ complexes $C_1,\dots,C_{5m}$ constructed as above are pairwise non-isomorphic. This follows from the fact that the graphs $\Gamma^\prime_1,\dots,\Gamma^\prime_{5m}$ are pairwise non-isomorphic, and each graph $\Gamma^\prime_{i}$ can be recovered, \textit{by construction}, from the data of the pairing maps which define $C_i$. Moreover, these complexes have trivial automorphism group. Indeed, let $\phi$ be an automorphism of $C_i$. Since $\phi$ preserves the number of self-identification of each simplex in $C_i$, it sends simplices with a single self-identification of facets to those with a single self-identification. 

In each complex $C_i$, we have two kinds of such simplices: the vertex simplices $S_v$, and those inside each edge sub-complex $E^\prime_k$. However, we notice that each simplex with a single self-identification in $E^\prime_k$ is paired with a simplex with two self-identifications. In contrast, each $S_v$ is paired with a simplex in some $E_k$ that has only one self-identification. Since $\phi$ preserves the glueing of the facets, each vertex sub-complex is mapped to a vertex sub-complex, and each edge sub-complex is mapped to an edge sub-complex. 

Thus, $\phi$ also acts as an automorphism of the asymmetric graph $\Gamma^\prime_i$. As a consequence of this fact, $\phi$ is required to fix all the vertex sub-complexes, as well as each copy $E^\prime_k$ of the edge sub-complex $E$. Since the edge sub-complexes have trivial automorphism group, the map $\phi$ has to be the identity on each of them. Since each vertex sub-complex has its free facets paired to some edge sub-complex $E^\prime_k$, $\phi$ is necessarily the identity on the vertex sub-complexes too.

Now we simply use the complexes $C_1,\dots,C_{5m}$ instead of the edge sub-complexes $E_k$ in the construction of Theorem \ref{teo:aut-triangulation} (recall that each complex $C_i$ has two free facets) in order to obtain a triangulation $\mathcal{T}$ with $\mathrm{Aut}\,\mathcal{T} \cong \mathrm{Aut}^{+}\,\mathcal{T} \cong G$. The number of simplices in such a triangulation can be computed using the following facts:
\begin{enumerate}
\item We have $v$ simplices coming from the modified Cayley graph $\Gamma^\prime$.
\item Each complex $C_i$ contains $8 N - 2$ simplices. Of these, $2 N$ simplices correspond to the vertex sub-complexes, and $6 N - 2$ simplices come from the edge complexes $E^\prime_1,\dots,E^\prime_{3 N - 1}$.
\end{enumerate}

The hyperbolic volume $V$ of the manifold $M_{\mathcal{T}}$ is therefore equal to 
\begin{equation}
v_\mathcal{B} \cdot (v + e\cdot (8 N - 2)),
\end{equation} 
where $v_\mathcal{B} = 4 \pi^2 / 3$ is the volume of the building block $\mathcal{B}$ constructed in Section \ref{sec:buildblock}. Therefore $V\sim N$ for $N\to \infty$. 

By choosing different collections of $5m$ pairwise non-isomorphic graphs $\Gamma_1, \dots, \Gamma_{5m}$ out of $f(N)$ totally available, we produce non-equivalent triangulations, and thus non-isometric manifolds. In fact, we have at least $\binom{f(N)}{5m}$ of them, and
\begin{align*}
\binom{f(N)}{5m} = \frac{f(N)!}{5m! \cdot (f(N)-5m)!} \sim  (f(N)-5m+1)\cdot (f(N)-5m+2)\cdot \dots \\ \ldots  \cdot (f(N)-2)\cdot (f(N)-1)\cdot f(N)
\sim f(N)^{5m} \sim C^{\, 5m\cdot N \cdot \log N}_0 \geq C^{\, N \log N}_0 \sim \\ \sim C^{\, V \log V}_1
\end{align*}
for some constants $C_0, C_1 > 1$, independent of the group $G$. Thus, $\rho_G(V) \geq C^{\, V \log {V}}$, for some $C>1$ and $V \geq V_0$, big enough. \prfend

\begin{rem}
By a classical result of Wang \cite{W}, in every dimension $n\geq4$ and for every $V>0$, there exist a finite number of hyperbolic $n$-manifolds of volume at most $V$, up to an isometry. Let $\rho_n(V)$ be the number of these manifolds. In \cite{BGLM} it is shown that, for every $n\geq4$, there exist constants $a,b\geq 0$ such that, for $V$ sufficiently large,

\begin{equation}
aV\log(V)\leq \log(\rho_n(V))\leq bV\log(V).
\end{equation}

By applying Theorem \ref{teo:superexponential}, it is easy to see that, in dimension four, an analogous bound holds for the functions $\rho_G(V)$: for any finite group $G$, there exist constants $a^{\prime}, b^{\prime}\geq 0$ independent of $G$ such that, for $V$ sufficiently large,
\begin{equation}
a^{\prime} V\log(V)\leq \log(\rho_G(V))\leq b^{\prime} V\log(V).
\end{equation}

\end{rem}

\newpage
 
\section{Appendix A}

Below we compute the volume of the rectified $5$-cell, realised as a four-dimensional ideal hyperbolic polytope.

By Heckman's formula \cite{H}, we have that
\begin{equation}
\mathrm{Vol}\,\mathcal{R} = \frac{1}{2}\,\cdot\mathrm{Vol}\,\mathbb{S}^4\,\cdot\chi_{orb}(\mathbb{H}^4/G),
\end{equation}
where $G = G(\mathcal{R})$ is the reflection group generated by the reflections in the facets of $\mathcal{R}$. Also, $\mathrm{Vol}\,\mathbb{S}^4 = 8\pi^2/3$.

Now we have to compute $\chi_{orb}(\mathbb{H}^4/G) = \chi(G)$, that is the Euler characteristic of the Coxeter group $G$. We suppose that $G$ is generated by the standard set of generators $S = \{s_1,\dots,s_m\}$, corresponding to the reflections in the facets of $\mathcal{R}$. Then, by a result of J.-P. Serre \cite{Serre},
\begin{equation}
\chi(G) = \sum_{T\subset S}\, \frac{(-1)^{|T|}}{|\langle T \rangle|},
\end{equation}
where we require $T$ to generate a finite special subgroup $\langle T \rangle \subset G$. 

In the reflection group $G = G(\mathcal{R})$, all possible non-trivial special finite subgroups are the stabilisers of its facets, two-dimensional faces, edges and vertices. Also, we have to count $\emptyset$ in the above formula, giving us a term of $1$.

We know that the rectified $5$-cell $\mathcal{R}$ has
\begin{enumerate}
\item 10 facets $P$, for each $Stab(P) \cong D_1 \cong \mathbb{Z}_2$ and $|D_1| = 2$.
\item 30 two-dimensional faces $F$. If $F$ is a face of a tetrahedral facet, then $Stab(F) \cong D_2$ (all the dihedral angles between tetrahedral and octahedral facets are right, tetrahedral facets do not intersect each other). If $F$ is a face of an octahedral facet, then $Stab(F)\cong D_3$ (octahedral facets intersect at the angle of $\pi/3$). Thus, we have 20 faces $F$ with $Stab(F)\cong D_2$ (all the faces of tetrahedral facets) and 30-20 = 10 remaining facets $F$ with $Stab(F)\cong D_3$. We know $|D_2|=4$ and $|D_3|=6$.
\item 30 edges $E$. Each edge $E$ is adjacent to one tetrahedral and two octahedral facets. Its stabiliser is the triangular reflection group $\Delta(2,2,3)$. We have $|\Delta(2,2,3)| = 12$. 
\item 10 vertices $V$. However, each vertex $V$ is ideal, and $Stab(V)$ is an affine reflection group. Thus, $Stab(V)$ is infinite and does not count.
\end{enumerate}

Finally, we compute
\begin{eqnarray}
\chi(G) = 1 - \frac{10}{|D_1|} + \left( \frac{20}{|D_2|} + \frac{10}{|D_3|} \right) - \frac{30}{|\Delta(2,2,3)|} =\\
\nonumber 1 - 10/2 + 20/4 + 10/6 - 30/12 = 1/6.
\end{eqnarray}
and, by Heckman's formula,
\begin{equation}
\mathrm{Vol}\,\mathcal{R} = \frac{1}{2}\,\cdot\frac{8\pi^2}{3}\,\cdot\frac{1}{6} = \frac{2\pi^2}{9}.
\end{equation}

\section{Appendix B}

In this appendix we provide the details of the SAGE routine computation \cite{K} used in the proof of Proposition~\ref{prop:cusp-length}. 

Let us recall that we need to show that the maximal cusp section $X$ of the manifold $M_{\mathcal{T}}$ produced by a triangulation $\mathcal{T}$ can be recovered from the length $h$ of the associated face cycle $c$ in $\mathcal{T}$. There are only three cases requiring computer aid. Namely, the cases when the cycle length equals $h=1$, and the respective return map $r_c$ is an even permutation. Then, by Lemma~\ref{lemma:cusptranslation}, the lifted action of the group $G$, the group of automorphisms of the Klein bottle $K$ preserving its tessellation by triangles, is induced by translation of the plane. Then, the translation lattice $L$ can be generated by $v_1 = (0, \sqrt{3}, 0)$, $v_2 = (3, 0, 0)$ and $v_3 = (1/2 + n, \sqrt{3}/2, 1)$, with $n=0,1,2$. We shall check that in either case there is no vector $v \in L$, such that $l(v) < l(v_2) = 3$, satisfying Conditions 1-3 from the proof of Proposition~\ref{prop:cusp-length}.

First of all, we find all possible vectors $v\in L$ of length $l(v) < 3$. Since every vector $v \in L$ can be represented as an integer linear combination of $v_1$, $v_2$ and $v_3$, we suppose
\begin{equation}
v = a\cdot v_1 + b\cdot v_2 + c\cdot v_3,
\end{equation}
and seek all possible $a,b,c \in \mathbb{Z}$ such that 
\begin{equation}
l^2(v) = w^T\, Q\, w \leq 9,
\end{equation}
where the matrix $Q$ is given by

\begin{equation}
\left(\begin{array}{ccc}
3 &0 &3/2\\
0 &9 &3n+3/2\\
3/2 &3n+3/2 &n^2+n+1/4
\end{array}\right)
\end{equation}

In each case $n=0,1,2$ we can compute the eigenvalues of $Q$ and conclude that $Q$ is positively defined. Thus, in some coordinate system $(x,y,z)$ it will define basically a new vector norm such that for each $v\in L$ we have
\begin{equation}
l^2(v) = \lambda_x\,x^2 + \lambda_y\,y^2 + \lambda_z\,z^2,
\end{equation}
where we suppose that $0 < \lambda_x \leq \lambda_y \leq \lambda_z$. 

Thus, it suffices to check all possible integer vectors $w$ in the ball of radius 
\begin{equation}
R = \floor*{\frac{3}{\sqrt{\lambda_x}}} + 1.
\end{equation}

The change of the coordinate system that brings the matrix $Q$ to its diagonal form is an orthogonal transformation, so the radius $R$ ball will be preserved under this coordinate change.

Thus, we first seek all possible integer solutions to the inequality
\begin{equation}
w^T\, Q\, w \leq 9, 
\end{equation}
with $w = (a, b, c) \in \mathbb{Z}^3$. 

Once we find such a vector $w$, we check if the quantity $h(w) = 3\sqrt{3}/l(w)^3$ is an integer for \textit{either} $w$ or $w/\sqrt{3}$. Indeed, we want to verify both possibilities: $v_1 = w$ or $v_2 = w$. 

The following SAGE routine performs all necessary computations.

\begin{verbatim}
n = 0;
Q = Matrix([[3,0,3/2],[0,9,3/2*(2*n+1)],[3/2,3/2*(2*n+1),7/4+(2*n+1)^2/4]]);
Q;
[  3   0 3/2]
[  0   9 3/2]
[3/2 3/2   2]
Q.eigenvalues();
[0.7345855820807065, 3.942467983256047, 9.32294643466325]
R = 3/sqrt(min(Q.eigenvalues()));
R;
3.500257982314587
R = floor(R) + 1;
R;
4
def lenQ(v): #returns the length of the vector (a,b,c) wrt 
             #the quadratic form defined by Q
    return sqrt((v*Q*v.column())[0]);
def is_int(x): #checks if number x is integer
    return (x - int(x) == 0);
def h(v): #computes the quantity h for vector v for face cycle of length 1
    return 3*sqrt(3)/(lenQ(v)**3);
def conditions(v): #verifies l(v) <= 3 and Condition 3 for vectors v 
                   #and v/sqrt(3)
    return (lenQ(v) <= 3) and (is_int(h(v)) or is_int(h(v/sqrt(3))));
l = list();
for a in range(-R,R):
    for b in range(-R,R):
        for c in range(-R,R):
            v = vector([a,b,c]);
            if lenQ(v)<>0 and conditions(v): l.append(v);
len(l);
4
l;
[(-1, 0, 0), (0, -1, 0), (0, 1, 0), (1, 0, 0)]    
\end{verbatim}

We get an analogous output in the cases $n=1$ and $n=2$. The coordinates $(a,b,c)$ of the output vectors mean that the only vectors in $L$ satisfying Conditions 1-3 are $\pm v_1$ and $\pm v_2$. In fact, we obtain some two dozens of solutions to the inequality $l(v)\leq 3$, and only Condition 3 sweeps out all but the desired ones.

\bigskip

\begin{flushleft}
\textit{Alexander Kolpakov\\
Department of Mathematics\\
40 St. George Street\\
Toronto ON\\
M5S 2E4 Canada\\}
\texttt{kolpakov dot alexander at gmail dot com}
\end{flushleft}

\medskip

\begin{flushleft}
\textit{Leone Slavich\\
Dipartimento di Matematica \\
Piazza di Porta San Donato 5\\
40126 Bologna, Italy\\}
\texttt{leone dot slavich at gmail dot com}
\end{flushleft}

\end{document}